\documentclass[MSNbibl,number,citesort,seceqn,dvips]{arxbj}

\aid{0}
\volume{20}
\issue{3}
\pubyear{2014}
\firstpage{1507}
\lastpage{1531}
\doi{10.3150/13-BEJ531} 

\makeatletter
\newcommand{\RMo}{\mathrm{o}}
\newcommand{\RMO}{\mathrm{O}}

\newcommand{\RMe}{\mathrm{e}}

\renewcommand{\mid}{|}

\newcommand{\mrmd}{\,\mathrm{d}}
\newcommand{\rrVert}{\Vert}
\newcommand{\rrvert}{\vert}
\newcommand{\llVert}{\Vert}
\newcommand{\llvert}{\vert}

\newtheorem{prop}{Proposition}[subsection]

\makeatother

\begin{document}
\begin{frontmatter}

\title{Asymptotic properties of adaptive maximum likelihood estimators
in latent variable models}
\runtitle{Adaptive ML estimators in GLLVM}

\begin{aug}
\author{\inits{S.}\fnms{Silvia} \snm{Bianconcini}\corref{}\ead[label=e1]{silvia.bianconcini@unibo.it}}
\runauthor{S. Bianconcini} 
\address{Department of Statistical Sciences, University of Bologna,
Via delle Belle Arti, 41-40126 Bologna, Italy.
\printead{e1}}
\end{aug}

\received{\smonth{10} \syear{2012}}
\revised{\smonth{5} \syear{2013}}

%
\begin{abstract}
Latent variable models have been widely applied in different fields of
research in which the constructs of interest are not directly
observable, so that one or more latent variables are required to reduce
the complexity of the data. In these cases, problems related to the
integration of the likelihood function of the model arise since
analytical solutions do not exist. In the recent literature, a
numerical technique that has been extensively applied to estimate
latent variable models is the adaptive Gauss--Hermite quadrature. It
provides a good approximation of the integral, and it is more
feasible than classical numerical techniques in presence of many latent
variables and/or random effects. In this paper, we formally
investigate the properties of maximum likelihood estimators based on
adaptive quadratures used to perform inference in generalized linear
latent variable models.
\end{abstract}

%

\end{frontmatter}

\section{Introduction}\label{sec1}

Models based on latent variables are used in many scientific fields,
particularly in social sciences. For instance, in psychology,
researchers often use concepts as intelligence and anxiety, that are
difficult to observe directly, but that can be indirectly measured by
surrogate data based on individual responses to a battery of tests. In
economics, welfare and poverty cannot be measured directly; hence
income, expenditure and various other indicators on households are used
as substitutes. Factor analysis is probably the best known latent
variable model, based on the assumption of multivariate normality for
the distribution of the manifest and latent variables. It has been
extended by numerous researchers in order to deal with survey data that
generally contain variables measured on binary, categorical or metric
scales, or combinations of the above.
Moustaki and Knott \cite{MouKno00} proposed a Generalized Linear
Latent Variable Model (GLLVM) framework that allows the distribution of
the manifest variables to belong to the exponential family, that
is either continuous or discrete variables. 


The purpose of GLLVM is to describe the relationship between a set of
responses or items $y_{1},\ldots, y_{p}$, and a set of latent
variables or factors $z_{1},\ldots, z_{q}$, that are fewer in number
than the observed variables. The factors are supposed to account for
the dependencies among the response variables in the sense that if the
factors are held fixed, then the observed variables are independent.
This is known as the assumption of conditional or local independence.
The conditional distribution of $y_{j}|\mathbf{z}$ ($\mathbf
{z}=[z_{1},\ldots, z_{q}]^{T}$) is taken from the exponential family (with
canonical link functions)
\[
g_{j}(y_{j}|\mathbf{z})=\exp\biggl\{
\frac{y_{j}(\alpha_{0j}+\bolds\alpha_{j}^{T} \mathbf
{z})-b_{j}(\alpha_{0j}+\bolds\alpha_{j}^{T} \mathbf{z})}{\phi
_{j}}+c_{j}(y_{j},\phi_{j}) \biggr\},\qquad
j=1,\ldots, p,
\]
where $\alpha_{0j}$ is the item-specific intercept, $\bolds\alpha
_{j}=[\alpha_{j1},\ldots,\alpha_{jq}]^{T}$ can be interpreted as factor
loadings of the model, and $\phi_{j}$ is the scale parameter, that is
of interest in the case of continuous observed components. 
The functions $b_{j}(\cdot)$ and $c_{j}(\cdot,\cdot)$ are known and
assume different forms according to the different nature of $y_{j}$.

Under the assumption of conditional independence, the joint marginal
distribution of the ma\-nifest variables is
%
%
\begin{equation}
\label{eq11} f(\mathbf{y}; \bolds\theta)=\int_{\mathbb{R}^q} g(
\mathbf{y}|\mathbf{z};\bolds\theta)h(\mathbf{z})\,\mathrm{d}\mathbf
{z}=\int
_{\mathbb{R}^q} \Biggl[\prod_{j=1}^{p}
g_{j}(y_{j}|\mathbf{z};\bolds\theta) \Biggr]h(
\mathbf{z})\,\mathrm{d}\mathbf{z}
\end{equation}
with $\mathbf{y}=[y_{1},\ldots,y_{p}]^{T}$, $\bolds\theta=[\alpha
_{01},\ldots, \alpha_{0p}, \bolds\alpha_{1}^{T},\ldots,
\bolds\alpha_{p}^{T}, \phi_{1},\ldots, \phi_{p}]^{T}$, and where
$\mathbf{z}$ is generally assumed to be multivariate standard normal,
but the independence assumption of the latent variables could be relaxed.

GLLVMs are designed as a flexible modelling approach. As a consequence,
they are rather complex models, and their statistical analysis presents
some difficulties due to the fact that the latent variables are not
observed. Maximum likelihood estimates in the GLLVM framework are
typically obtained by using standard maximization algorithms, such as
the EM and the Newton--Raphson algorithms
(Moustaki and Knott \cite{MouKno00}, Huber, Ronchetti and
Victoria-Feser \cite
{HubRonFe04}).
In both cases, the latent variables must be integrated out from the
likelihood function, and numerical techniques have to be applied.
Moustaki and Knott \cite{MouKno00} proposed the use of the
Gauss--Hermite (GH) quadrature
as a numerical approximation method. Although this is feasible in
fairly simple models and tends to work well with moderate sample sizes,
its application is often unfeasible when the number of latent variables
increases. Moreover, GH can completely miss the maximum for certain
functions and can be inefficient in other cases. To overcome these
limitations, the Adaptive Gauss--Hermite (AGH) quadrature has become
very popular in the latent variable literature. It allows to get a
better approximation of the integral by adjusting the quadrature
locations with specific features of the posterior density of the latent
variables given the observations. Developed in the Bayesian context by
Naylor and Smith \cite{NaySmi82}, it has been extended by several
authors to deal with
generalized linear mixed models. In particular,
Schilling and Bock \cite{SchBo05} applied the AGH quadrature to
approximate marginal likelihoods in IRT models with binary data,
whereas Rabe-Hesketh, Skrondal and Pickles \cite{HesSkPi05} analyzed
its behavior
for generalized linear latent and mixed models. Furthermore,\vadjust{\goodbreak} 
Joe \cite{Joe08} compared the AGH with the Laplace
approximation for a variety of discrete response mixed models. He found
that the Laplace approximation becomes less adequate as the degree of
discreteness increases and suggests using AGH with binary and ordinal
data. On this regard, we recall that several approaches have been
proposed to overcome the main limitations of the Laplace approximation.
In latent Gaussian models
(Rue and Held \cite{Rue05}), the Integrated Nested Laplace
Approximation (INLA) has become very popular to perform Bayesian
inference with non-Gaussian observations
(Rue, Martino and Chopin~\cite{Rue09}). This procedure combines Laplace
approximations with numerical integration to provide a fast and
accurate method for approximating the predictive density of the latent
variables/random effects. It is also a valuable tool in practice via
the R-package \emph{R-INLA}
(Martins \textit{et al.} \cite{Rue13}).

The adaptive Gauss--Hermite quadrature is implemented in many
statistical software used to fit GLLVM, such as in the function \emph
{gllamm} in STATA
(Rabe-Hesketh and Skrondal \cite{RaSk12}), in MPLUS (Muthen and Muthen
\cite{MuMu10}), and in the PROC NLMIXED in SAS
(Lesaffre and Spiessens \cite{LeSp01}). However, to the best of our
knowledge, inferential issues on the properties of the estimators based
on the adaptive quadrature have not been addressed in the literature.
In this paper, we formally investigate these theoretical properties as
function of both the sample size and the number of observed variables.
Our results generalize those by Huber, Ronchetti and Victoria-Feser
\cite
{HubRonFe04}, who analyzed the properties of classical Laplace-based
estimators in GLLVM. Indeed, we show that the adaptive Gauss--Hermite
quadratures share the same error rate of the higher (than one) order
Laplace approximation.


The paper is organized as follows. In Section~\ref{sec2}, we discuss the
estimation of GLLVMs when the adaptive Gauss--Hermite quadrature is
applied to approximate integrals. In Section~\ref{sec3}, the relationship
between AGH quadratures and the Laplace approximation is analyzed, and
the asymptotic properties of the adaptive Maximum Likelihood (ML) estimators are derived. A
simulation study is implemented in Section~\ref{sec4} to analyze the finite
sample properties of the estimators. Finally, in Section~\ref{sec5}, a brief
summary on the main findings of the paper is provided.\vspace*{-3pt}

\section{Estimation based on adaptive Gauss--Hermite
quadrature}\label{sec2}\vspace*{-3pt}
Maximum Likelihood (ML) estimates in the GLLVM framework are typically
obtained by using either the EM or the Newton--Raphson algorithms. The
key component for applying both the algorithms is the score vector of
the observed data log-likelihood function.
For a random sample of size $n$, the latter is defined as
%
%
\begin{eqnarray}
\label{loglike}
\ell(\bolds\theta)&=&\sum
_{l=1}^{n}\log f(\mathbf{y}_{l};
\bolds\theta)
\nonumber\\[-2pt]
&=&\sum_{l=1}^{n} \log\int
_{\mathbb{R}^{q}} \prod_{j=1}^{p}\exp
\biggl[\frac
{y_{jl} (\alpha_{0j}+\bolds\alpha_{j}^{T} \mathbf
{z}_{l})-b_{j}(\alpha_{0j}+\bolds\alpha_{j}^{T} \mathbf
{z}_{l})}{\phi_{j}}+c_{j}(y_{jl},\phi_{j})
\biggr]
\\[-2pt]
&&\hspace*{58.5pt}{}\times(2\pi)^{-q/2} \exp\biggl[-\frac{1}{2}\mathbf{z}_{l}^{T}
\mathbf{z}_{l} \biggr]\mrmd \mathbf{z}_{l}.\nonumber
\end{eqnarray}
It is easily shown that the score vector corresponding to expression
(\ref{loglike}) equals
%
%
\begin{eqnarray}\label{score}
S(\bolds\theta)&=&\frac{\partial\ell(\bolds\theta)}{\partial
\bolds\theta}=\sum
_{l=1}^{n}\frac{\partial}{\partial
\bolds\theta} \log f(
\mathbf{y}_{l};\bolds\theta)
\nonumber\\
&=&\sum_{l=1}^{n}\frac{1}{f(\mathbf{y}_{l}; \bolds\theta)}
\int_{\mathbb{R}^{q}}\frac{\partial}{\partial\bolds\theta} \bigl
[g(\mathbf{y}_{l}
\mid\mathbf{z}_{l}; \bolds\theta)h(\mathbf{z}_{l})
\bigr]\mrmd \mathbf{z}_{l}
\\
&=& \sum_{l=1}^{n}
\frac{\int_{\mathbb{R}^{q}}S_{l}(\bolds\theta;
\mathbf{z}_{l})g(\mathbf{y}_{l} \mid\mathbf{z}_{l}; \bolds\theta
)h(\mathbf{z}_{l})\mrmd \mathbf{z}_{l}}{\int_{\mathbb{R}^{q}}g(\mathbf
{y}_{l} \mid\mathbf{z}_{l}; \bolds\theta)h(\mathbf
{z}_{l})\mrmd \mathbf{z}_{l}}
\nonumber\\
&=&\sum_{l=1}^{n}\int
_{\mathbb{R}^{q}}S_{l}(\bolds\theta; \mathbf
{z}_{l})h(\mathbf{z}_{l} \mid\mathbf{y}_{l};
\bolds\theta)\mrmd \mathbf{z}_{l}=\sum_{l=1}^{n}
E_{\mathbf{z}|\mathbf
{y}}\bigl[S_{l}(\bolds\theta; \mathbf{z}_{l})
\bigr],\nonumber
\end{eqnarray}
where $S_{l}(\bolds\theta; \mathbf{z}_{l})$ denotes the
complete-data score vector given by $\partial\log f(\mathbf{y}_{l},
\mathbf{z}_{l}; \bolds\theta)/\partial\bolds\theta
=\partial[\log g(\mathbf{y}_{l} \mid\mathbf{z}_{l}; \bolds
\theta) + \log h(\mathbf{z}_{l}) ]/\partial\bolds\theta$.
In words, the observed data score vector is expressed as the expected
value of the complete-data vector with respect to $h(\mathbf{z}_{l}
\mid\mathbf{y}_{l}; \bolds\theta)$, that is the posterior
distribution of the latent variables given the observations. This
implies that (\ref{score}) plays a double role. If the score
equations are solved with respect to $\bolds\theta$, with
$h(\mathbf{z}_{l} \mid\mathbf{y}_{l}; \bolds\theta)$ fixed at
the $\bolds\theta$-value of the previous iteration, then this
corresponds to the EM algorithm, whereas, if the score equations are
solved with respect to $\bolds\theta$ considering $h(\mathbf
{z}_{l} \mid\mathbf{y}_{l}; \bolds\theta)$ also as a function of
$\bolds\theta$, then this corresponds to a direct maximization of
the observed data log-likelihood $\ell(\bolds\theta)$. As we
shall discuss further, based on this appealing feature, the estimators
derived by applying either of these two algorithms will share the same
theoretical properties.

Equation (\ref{score}) involves ratios of multidimensional integrals which
cannot be solved analytically, except when all the $g_{j}(y_{jl}|\mathbf
{z}_{l}; \bolds\theta)$ are normal. Consequently, an
approximation of these integrals is needed, on which the bias and
variance of resulting estimators will depend.
In this paper, we study the properties of ML estimators based on the
adaptive Gauss--Hermite approximation of integrals. This technique
consists of adjusting the quadrature locations with specific features
of the posterior density of the latent variables given the
observations. This provides a better approximation of the function to
be integrated.
Naylor and Smith \cite{NaySmi82} took the mean vector and
covariance matrix of the normal density approximating the integrand to
be the posterior mean and covariance matrix. Unfortunately, these
posterior moments are not known exactly, but must themselves be
obtained using adaptive quadratures. Integration is therefore
iterative. To overcome this limitation, Liu and Pierce \cite
{LiuPie94} proposed an alternative procedure that consists in computing
the mode of the integrand and its curvature (inverse of the Hessian
matrix) at the mode, so that numerical integration is avoided. In this
case, the adaptive quadrature, when applied using one abscissa, is
equivalent to the classical Laplace approximation, and its behavior has
been analyzed in several papers on generalized linear models
(Pinhero and Bates \cite{PiBa95}, Schilling and Bock \cite{SchBo05},
Skrondal and Rabe-Hesketh \cite{SkRa04}, Joe \cite{Joe08}).

The application of the adaptive quadrature requires to rewrite (\ref
{eq11}) as follows
%
%
\begin{equation}
\label{AGH} f(\mathbf{y}_{l}; \bolds\theta)=\int
_{\mathbb{R}^{q}} \frac{g(\mathbf{y}_{l} \mid\mathbf{z}_{l}; \bolds
\theta)h(\mathbf{z}_{l})}{h_{1}(\mathbf{z}_{l};
\hat{\mathbf{z}}_{l}, \bolds
\Psi_{l})}h_{1}(\mathbf{z}_{l};
\hat{\mathbf{z}}_{l}, \bolds\Psi_{l})\mrmd
\mathbf{z}_{l},
\end{equation}
where
$h_{1}(\cdot;\hat{\mathbf{z}}_{l}, \bolds\Psi_{l})$
is a multivariate normal density with first and second moments
%
%
\begin{eqnarray}
\label{mode} \hat{\mathbf{z}}_{l}&=&\arg\max_{\mathbf{z}_{l}\in
\mathbb{R}^{q}}
\bigl[\log g(\mathbf{y}_{l}|\mathbf{z}_{l};\bolds
\theta)+ \log h(\mathbf{z}_{l}) \bigr],
\\
%
%
\label{hessian} \bolds\Psi_{l} &=& \biggl(- \frac{\partial^{2} [\log
g(\mathbf{y}_{l} \mid\mathbf{z}_{l};
\bolds\theta)+ \log h(\mathbf{z}_{l}) ]}{\partial
\mathbf{z}_{l}^{T}\,\partial\mathbf{z}_{l}}
\biggr)\bigg|_{\mathbf{z}_{l}
=\hat{\mathbf{z}}_{l}}^{-1}.
\end{eqnarray}
A cartesian product rule based on the classical Gauss--Hermite
quadrature is then applied so that the integrals have to be defined
with respect to uncorrelated variables $\tilde{\mathbf{z}}_{l}$. Based
on the Cholesky factorization of the covariance matrix $\bolds
\Psi_{l} = \mathbf{T}_{l}\mathbf{T}_{l}^{T}$, expression (\ref{AGH})
can be rewritten as
\[
f(\mathbf{y}_{l}; \bolds\theta)=2^{{q}/{2}}|
\mathbf{T}_{l}|\int_{\mathbb{R}^{q}} g(\mathbf{y}_{l}|
\sqrt{2}\mathbf{T}_{l}\tilde{\mathbf{z}}_{l} + \hat{
\mathbf{z}}_{l}; \bolds\theta)h(\sqrt{2}\mathbf{T}_{l}
\tilde{\mathbf{z}}_{l} + \hat{\mathbf{z}}_{l})\exp\bigl[
\tilde{\mathbf{z}}_{l}^{T}\tilde{\mathbf{z}}_{l}
\bigr]\exp\bigl[-\tilde{\mathbf{z}}_{l}^{T}\tilde{\mathbf
{z}}_{l} \bigr]\mrmd \tilde{\mathbf{z}}_{l},
\]
such that the AGH approximation of the density $f(\mathbf{y}_{l};
\bolds\theta), l=1,\ldots, n$, is given by
%
%
\begin{equation}
\label{marginal} \tilde{f}(\mathbf{y}_{l}; \bolds
\theta)=2^{q/2}|\mathbf{T}_{l}| \sum
_{t_{1},\ldots, t_{q}} g\bigl(\mathbf{y}_{l}|\mathbf
{z}^{*}_{l, t_{1},\ldots, t_{q}}; \bolds\theta\bigr)h\bigl(\mathbf
{z}^{*}_{l, t_{1},\ldots, t_{q}}\bigr)w_{t_{1}}^{*}\cdots
w_{t_{q}}^{*},
\end{equation}
where $\sum_{t_{1},\ldots, t_{q}}=\sum_{t_{1}=1}^{k} \cdots\sum
_{t_{q}=1}^{k}$, being $k$ the number of quadrature points selected for
each latent variable, $\mathbf{z}^{*}_{l, t_{1},\ldots, t_{q}}=(z_{l,
t_{1}}^{*},\ldots, z_{l, t_{q}}^{*})^{T}= \sqrt{2}\mathbf
{T}_{l}(z_{t_{1}},\ldots, z_{t_{q}})^{T}+ \hat{\mathbf{z}}_{l}$ and
$w_{t_{k}}^{*}=w_{t_{k}}\exp[z_{t_{k}}^{2}]$ are the AGH nodes and
weights, respectively, with $z_{t_{k}}$ being the classical GH nodes
and $w_{t_{k}}$, $k=1,\ldots,q$, the corresponding weights.

From (\ref{marginal}), we obtain the approximated log-likelihood function
%
%
\begin{eqnarray}
\label{approlike}
\tilde{\ell}(\bolds\theta)&=& \sum
_{l=1}^{n} \log\Biggl[2^{q/2}|
\mathbf{T}_{l}| \sum_{t_{1},\ldots, t_{q}} \prod
_{j=1}^{p} \exp\biggl(\frac{y_{jl} (\alpha_{0j}+\bolds\alpha
_{j}^{T}\mathbf{z}_{l, t_{1},\ldots,t_{q}}^{*} )-b_{j} (\alpha
_{0j}+\bolds\alpha_{j}^{T}\mathbf{z}_{l, t_{1},\ldots,t_{q}}^{*}
)}{\phi_{j}}\nonumber\\
&&\hspace*{128.5pt}{}+c_{j}(y_{jl},
\phi_{j}) \biggr)
\\
&&\hspace*{90pt}{}\times(2\pi)^{-q/2}\exp\biggl(-\frac{1}{2}\mathbf{z}_{l, t_{1},\ldots,t_{q}}^{*T}
\mathbf{z}_{l, t_{1},\ldots,t_{q}}^{*} \biggr)w_{t_{1}}^{*}
\cdots w_{t_{q}}^{*}\Biggr].\nonumber
\end{eqnarray}
%
The estimators of the model parameters are found by equating the
corresponding derivatives of (\ref{approlike}) to zero, that is
%
%
\begin{eqnarray}
\label{sco1}
\tilde{S}(\bolds\theta)&=&\frac
{\partial\tilde{\ell}(\bolds\theta)}{\partial\bolds\theta
}=\sum
_{l=1}^{n}\frac{1}{\tilde{f}(\mathbf{y}_{l}; \bolds\theta
)}\,
\frac{\partial\tilde{f}(\mathbf{y}_{l}; \bolds\theta
)}{\partial\bolds\theta}
\nonumber\\
&=&\sum_{l=1}^{n}\frac{\sum_{t_{1},\ldots, t_{q}}S_{l}(\bolds
\theta; \mathbf{z}^{*}_{l,
t_{1},\ldots, t_{q}}) g(\mathbf{y}_{l}|\mathbf{z}^{*}_{l, t_{1},\ldots,
t_{q}}; \bolds\theta)h(\mathbf{z}^{*}_{l, t_{1},\ldots,
t_{q}})w_{t_{1}}^{*}\cdots w_{t_{q}}^{*}}{\sum_{t_{1},\ldots, t_{q}}
g(\mathbf{y}_{l}|\mathbf{z}^{*}_{l, t_{1},\ldots, t_{q}}; \bolds
\theta)h(\mathbf{z}^{*}_{t_{1},\ldots, t_{q}})w_{t_{1}}^{*}\cdots
w_{t_{q}}^{*}}
\\
&=&\sum_{l=1}^{n}\tilde{E}_{\mathbf{z}| \mathbf{y}}
\bigl[S_{l}(\bolds\theta; \mathbf{z}_{l}) \bigr]=0,\nonumber
\end{eqnarray}
where, specifically,
\begin{eqnarray*}
S_{l}\bigl(\alpha_{0j}; \mathbf{z}^{*}_{l, t_{1},\ldots, t_{q}}
\bigr)&=& \frac
{1}{\phi_{j}} \biggl[y_{jl}-\frac{\partial b_{j} (\alpha
_{0j}+\bolds\alpha_{j}^{T}\mathbf{z}_{l, t_{1},\ldots,t_{q}}^{*}
)}{\partial\alpha_{0j}} \biggr],
\\
S_{l}\bigl(\bolds\alpha_{j}; \mathbf{z}^{*}_{l, t_{1},\ldots,
t_{q}}
\bigr)&=& \frac{\mathbf{z}^{*}_{l, t_{1},\ldots, t_{q}}}{\phi_{j}}
\biggl
[y_{jl}-\frac{\partial b_{j} (\alpha_{0j}+\bolds\alpha
_{j}^{T}\mathbf{z}_{l, t_{1},\ldots,t_{q}}^{*} )}{\partial
\bolds\alpha_{j}} \biggr]
\end{eqnarray*}
and
\[
S_{l}\bigl(\phi_{j}; \mathbf{z}^{*}_{l, t_{1},\ldots, t_{q}}
\bigr)= -\frac
{1}{\phi_{j}^{2}} \bigl[y_{jl} \bigl(\alpha_{0j}+
\bolds\alpha_{j}^{T}\mathbf{z}_{l, t_{1},\ldots,t_{q}}^{*}
\bigr) - b_{j} \bigl(\alpha_{0j}-\bolds
\alpha_{j}^{T}\mathbf{z}_{l, t_{1},\ldots,t_{q}}^{*} \bigr)
\bigr]+\frac{\partial c_{j}(y_{jl}, \phi
_{j})}{\partial\phi_{j}}.
\]
Equations (\ref{sco1}) provide a set of estimating equations defining
the estimators for the model parameters. 
The same equations are derived in the $E$-step of the EM algorithm, in
which the AGH quadrature is applied to approximate the $E$-step
expectations (\ref{score}). In the $M$-step, as in the direct
maximization algorithm, improved estimates for the model parameters are
obtained by maximizing the approximated expected score functions (\ref
{sco1}). For the scale parameter $\phi_{j}$, closed form expressions
can be derived, whereas, for the other parameters, a Newton Raphson
iterative scheme is used in order to solve the corresponding nonlinear
maximum likelihood equations. In the derivation of the estimating
equations, the model has been kept as general as possible without
specifying the conditional distributions $g_{j}(y_{j}|\mathbf{z};
\bolds\theta)$. In Appendix \ref{app3}, we give specific expressions for
the quantity that are used in the log-likelihood function (\ref
{loglike}) and in the score functions (\ref{sco1}) for binary manifest
variables, whereas we refer to Bianconcini and Cagnone \cite{BiaCag12}
for count and categorical observed variables.

\section{Statistical properties of the AGH-based estimators}\label{sec3}
To investigate the asymptotic properties of the maximum likelihood
estimators based on the adaptive Gauss--Hermite quadrature, the error
rate associated to the approximation (\ref{sco1}) has to be
determined. Liu and Pierce \cite{LiuPie94} analyzed the
asymptotic behavior of the AGH
when it is used to approximate unidimensional integrals. Based on the
fact that when applied with only one node it results in the Laplace
approximation to the integral
(de Bruijn \cite{DeBr81}, Barndorff-Nielsen and Cox~\cite{BarCox89}),
they proved that
the adaptive quadrature based on $k$ points can be alternatively
thought as a higher (than one) order Laplace approximation. We now
generalize this result to the multidimensional integral (\ref{eq11})
as well as to the ratio of integrals (\ref{score}), and we analyze the
asymptotic accuracy of the corresponding Laplace approximations. The
behavior of the latter for multidimensional integrals was studied by
Barndorff-Nielsen and Cox \cite{BarCox89}, Shun and McCullagh \cite
{ShuMc95}, Shun \cite{Shu97}, and recently by
Evangelou, Zhu and Smith \cite{Eva11} for spatial generalized linear
mixed models.
Similarly, Raudenbush, Yang and Yosef \cite{Rade00} considered
improvements of
the standard Laplace approximation obtained by incorporating higher
order derivatives of the integrand.

For the derivations illustrated here, we follow the notation of
Shun and McCullagh \cite{ShuMc95} based on summation convention.
Hence, an index that appears as a subscript and as a superscript
implies a summation over all possible values of that index. We will
denote the components of a vector sometimes by subscripts and sometimes
by superscripts. 
The ($i,j$)th component of a matrix $\mathbf{A}$ will be written as
$a_{ij}$ and its inverse (when exists) will have components $a^{ij}$.
For any real function $f(\mathbf{z}), \mathbf{z} \in\mathbb{R}^{q}$,
its derivative with respect to the $i$th component of $\mathbf{z}$ is
denoted by a subscript, that is, $f_{i}(\mathbf{z})=\frac
{\partial f(\mathbf{z})}{\partial z_{i}}$, $f_{ij}(\mathbf{z})=\frac
{\partial^{2} f(\mathbf{z})}{\partial z_{i}\,\partial z_{j}}$, and, more
generally, $f_{i_{1},\ldots, i_{2m}}(\mathbf{z})=\frac{\partial
^{2m}f(\mathbf{z})}{\partial z_{i_{1}}\cdots\partial z_{i_{2m}}}$. In
order to keep the notation as light as possible, we omit the individual
subscript $l$.

\subsection{Relationship with the Laplace approximation}\label{sec3.1}

The AGH quadrature implemented here is based on a tensor product of $q$
univariate Gaussian quadratures based on the same number of quadrature
points. In each dimension, the approximation (\ref{marginal}) is exact
for polynomials of degree $2k+1$ or less. Hence, it provides a good
approximation of the integral (\ref{eq11}) if the ratio $\nu(\mathbf
{z})=\frac{g(\mathbf{y}|\mathbf{z}; \bolds\theta)h(\mathbf
{z})}{h_{1}(\mathbf{z};\hat{\mathbf{z}}, \bolds\Psi)}$ can be
approximated well by a $q$-variate polynomial, where the maximum
exponent of all the monomials is at most $2k+1$
(Tauchen and Hussey \cite{TauHus91}), in the region where the
integrand is substantial. It follows that the effectiveness of the
adaptive Gauss--Hermite approximation (\ref{marginal}) can be evaluated
by considering the Taylor series expansion of $\nu(\mathbf{z})$ around
the mode $\hat{\mathbf{z}}$, that is, 
%
%
\begin{equation}
\label{taylor} \nu(\mathbf{z})=\nu(\hat{\mathbf{z}}) \Biggl[1+\sum
_{m=3}^{\infty} \frac
{1}{m!}c_{i_{1},\ldots,i_{m}}(\hat{
\mathbf{z}}) (\mathbf{z}-\hat{\mathbf{z}})^{i_{1},\ldots,i_{m}}
\Biggr],
\end{equation}
where $(i_{1},\ldots, i_{m})$ is a set of $m$ indices, $c_{i_{1},\ldots
,i_{m}}(\hat{\mathbf{z}})=\frac{\nu_{i_{1},\ldots,i_{m}}(\hat{\mathbf
{z}})}{\nu(\hat{\mathbf{z}})}$, $\nu_{i_{1},\ldots,i_{m}}(\hat{\mathbf
{z}})$ denotes the partial derivatives of order $m$ of $\nu$ with
respect to $z_{i_{1}},\ldots, z_{i_{m}}$ evaluated at the mode $\hat
{\mathbf{z}}$, whereas $(\mathbf{z}-\hat{\mathbf{z}})^{i_{1},\ldots
,i_{m}}$ refers to specific components of the vector $(\mathbf{z}-\hat
{\mathbf{z}})$. The coefficients $c_{i_{1}}$ and $c_{i_{1},i_{2}}$ are
zero due to the choice of $h_{1}(\cdot; \hat{\mathbf{z}}, \bolds
\Psi)$.

Substituting the expansion (\ref{taylor}) into the integral
(\ref{eq11}), we obtain the exact solution
%
%
\begin{equation}
\label{taylor2} f(\mathbf{y}; \bolds\theta)=\nu(\hat{\mathbf{z}})
\Biggl[1+
\sum_{m=2}^{\infty} \sum
_{Q} \frac{1}{(2m)!}c_{i_{1},\ldots,i_{2m}}(\hat{\mathbf{z}})
\bolds\nu^{q_{1}}(\hat{\mathbf{z}})\cdots\bolds
\nu^{q_{m}}(\hat{\mathbf{z}}) \Biggr],
\end{equation}
where the second sum is over the partition $Q=q_{1}|\cdots|q_{m}$ of
$2m$ indices into $m$ blocks, each of size 2, and $\bolds\nu
^{q_{k}}(\hat{\mathbf{z}}), k=1,\ldots,m$, are components of the
covariance matrix $\bolds\Psi$. The Gauss--Hermite quadrature, for
which $k$ quadrature points are selected for each dimension, would be
exact if the partial derivatives beyond the $2(k+1)$ order in (\ref
{taylor2}) are zero, that is,
%
%
\begin{equation}
\label{taylor3} f(\mathbf{y}; \bolds\theta)=\nu(\hat{\mathbf{z}})
\Biggl[1+
\sum_{m=2}^{k} \sum
_{Q} \frac{1}{(2m)!}c_{i_{1},\ldots,i_{2m}}(\hat{\mathbf{z}})
\bolds\nu^{q_{1}}(\hat{\mathbf{z}})\cdots\bolds\nu
^{q_{m}}(\hat{\mathbf{z}}) \Biggr].
\end{equation}
To determine the asymptotic order of the approximation (\ref{taylor3}),
its relationship with the higher order Laplace approximation of
multidimensional integrals has to be taken into account. At this
regard, the integral (\ref{eq11}) has to be rewritten as
%
%
\begin{equation}
\label{FLA}f(\mathbf{y}; \bolds\theta)=\int_{\mathbb{R}_{q}}
\RMe ^{ [-L(\mathbf{z}) ]}\mrmd \mathbf{z},
\end{equation}
where $L(\mathbf{z})=- [\log g(\mathbf{y}| \mathbf{z}; \bolds
\theta) + \log h(\mathbf{z}) ]$, such that $L(\mathbf{z})=\RMO(p)$.
Assuming that $L(\mathbf{z})$ has a unique minimum $\hat{\mathbf{z}}$,
Shun and McCullagh \cite{ShuMc95} suggested the following expansion
around that minimum
\[
L(\mathbf{z})=L(\hat{\mathbf{z}})+\sum_{m=2}^{\infty}
\frac
{1}{m!}L_{i_{1},\ldots,i_{m}}(\hat{\mathbf{z}}) (\mathbf{z}-\hat{\mathbf
{z}})^{i_{1},\ldots, i_{m}}
\]
and applying the exponential function
\[
\RMe^{-L(\mathbf{z})}=(2\pi)^{q/2}|\bolds\Psi|^{1/2}\RMe^{-L(\hat
{\mathbf{z}})}h_{1}(
\mathbf{z}; \hat{\mathbf{z}}, \bolds\Psi)\exp\Biggl[\sum
_{m=3}^{\infty}\frac{(-1)}{m!}L_{i_{1},\ldots,i_{m}}(\hat{
\mathbf{z}}) (\mathbf{z}-\hat{\mathbf{z}})^{i_{1},\ldots, i_{m}}
\Biggr],
\]
where $h_{1}(\cdot;\hat{\mathbf{z}}, \bolds\Psi)$ is a
multivariate normal density with moments given in (\ref{mode}) and
(\ref{hessian}). Based on exlog relations, the higher order term can be
expressed as follows
\begin{eqnarray*}
&&
\exp\Biggl[\sum_{m=3}^{\infty}
\frac{(-1)}{m!}L_{i_{1},\ldots,i_{m}}(\hat{\mathbf{z}}) (\mathbf
{z}-\hat{
\mathbf{z}})^{i_{1},\ldots, i_{m}} \Biggr]\\
&&\quad=1 - \sum_{m=3}^{\infty}
\sum_{P}\frac{(-1)^{t}}{m!}L_{p_{1}}(\hat{
\mathbf{z}}) \cdots L_{p_{t}}(\hat{\mathbf{z}}) (\mathbf{z}-\hat{\mathbf
{z}})^{i_{1},\ldots, i_{m}}
\end{eqnarray*}
such that the exact solution of the integral (\ref{FLA}) is given by
%
%
\begin{equation}
\label{Ei2} (2\pi)^{q/2}|\bolds\Psi|^{
{1/2}}\RMe^{-L(\hat{\mathbf
{z}})}
\Biggl[1-\sum_{m=2}^{\infty}\sum
_{P,Q}\frac
{(-1)^{t}}{(2m)!}L_{p_{1}}(\hat{\mathbf{z}})
\cdots L_{p_{t}}(\hat{\mathbf{z}})L^{q_{1}}(\hat{\mathbf{z}})
\cdots L^{q_{m}}(\hat{\mathbf{z}}) \Biggr],
\end{equation}
where the second sum is over all partitions $P,Q$, such that
$P=p_{1}|\cdots|p_{t}$ is a partition of $2m$ indices into $t$ blocks,
each of size 3 or more, and $Q=q_{1}|\cdots|q_{m}$ is a partition of
$2m$ indices into $m$ blocks, each of size 2. Each component
$L^{q_{k}}(\hat{\mathbf{z}}), k=1,\ldots, m$, refers to specific
elements of $\bolds\Psi$. As shown in the Appendix \ref{app1}, the exact
solution (\ref{Ei2}) is equivalent to the one derived in (\ref
{taylor2}). It follows that the asymptotic order of the AGH
approximation can be derived by truncating at $m=k$ the expansion (\ref
{Ei2}), and by analyzing the asymptotic order associated to the
bipartition ($P,Q$) related to $m=k+1$. For fixed $q$, the usual
asymptotic order of the term corresponding to the bipartition $(P,Q)$
in (\ref{Ei2}) is $\RMO(p^{t-m})$. It follows that the asymptotic error of
the AGH based on $k$ quadrature points is the same associated to the
bipartition $(P,Q)$ of $2(k+1)$ indices, that is, $\RMO (p^{- [
{k}/{3}+1 ]} )$ (see in Appendix \ref{app1} for more details).

It has to be noticed that when AGH quadratures are applied in the
estimation of GLLVM, we need to approximate ratios of integrals as
shown in (\ref{score}). The fully exponential solution (\ref{Ei2})
cannot be applied to the integral at the numerator, since the score
functions $S(\bolds\theta; \mathbf{z})$ are not necessarily
positive. The integral has to be written in the standard form
(Tierney, Kass and Kadane \cite{TieKa89}, Evangelou, Zhu and Smith
\cite{Eva11})
\[
\int_{\mathbb{R}^{q}}\RMe^{-L(\mathbf{z})}S(\bolds\theta; \mathbf{z})\mrmd
\mathbf{z}.
\]
Beyond the Taylor series expansion of $L(\mathbf{z})$ around its
minimum $\hat{\mathbf{z}}$, we have to consider a similar expansion of
$S$ around the same point, that is,
\[
S(\bolds\theta; \mathbf{z})=\sum_{m=0}^{\infty}S_{j_{1},\ldots,j_{m}}(
\bolds\theta; \hat{\mathbf{z}}) (\mathbf{z} - \hat{\mathbf
{z}})^{j_{1},\ldots, j_{m}}.
\]
Following Evangelou, Zhu and Smith \cite{Eva11}, it can be shown that %
\begin{eqnarray*}
&&
\int_{\mathbb{R}^{q}}S(\bolds\theta; \mathbf{z})g(
\mathbf{y}|\mathbf{z}; \bolds\theta)h(\mathbf{z})\mrmd \mathbf{z}\\
&&\quad=(2
\pi)^{q/2}|\bolds\Psi|^{
{1/2}}\RMe^{-[L(\hat{\mathbf{z}})]}
\\
&&\qquad{}\times\Biggl[\sum_{m=0}^{\infty} \sum
_{s=0}^{2m}\sum
_{P,Q}\frac{(-1)^{t}}{(2m)!}S_{j_{1},\ldots,j_{s}}(\bolds\theta;
\hat{\mathbf{z}})L_{p_{1}}(\hat{\mathbf{z}}) \cdots L_{p_{t}}(
\hat{\mathbf{z}})L^{q_{1}}(\hat{\mathbf{z}})\cdots L^{q_{m}}(\hat{
\mathbf{z}})\Biggr],
\end{eqnarray*}
where $P$ is a partition of $2m-s$ indices into $t$ blocks, each of
size 3 or more, and $Q$ is a partition of the same indices together
with $\{j_{1},\ldots,j_{s}\}$ into $m$ blocks of size 2. Note that $P$
and $Q$ do not need to be connected. It follows that the exact Laplace
solution of the expected score function (\ref{score}) results%
%
%
\begin{equation}
\frac{\sum_{m=0}^{\infty} \sum_{s=0}^{2m}\sum_{P,Q}
({(-1)^{t}}/({2m!}))S_{j_{1},\ldots,j_{s}}(\bolds\theta; \hat{\mathbf
{z}})L_{p_{1}}(\hat{\mathbf{z}}) \cdots L_{p_{t}}(\hat{\mathbf
{z}})L^{q_{1}}(\hat{\mathbf{z}})\cdots L^{q_{m}}(\hat{\mathbf{z}})}{\sum
_{m=0}^{\infty}\sum_{P,Q}({(-1)^{t}}/({2m!}))L_{p_{1}}(\hat{\mathbf
{z}}) \cdots L_{p_{t}}(\hat{\mathbf{z}})L^{q_{1}}(\hat{\mathbf
{z}})\cdots L^{q_{m}}(\hat{\mathbf{z}})}.
\end{equation}
It will be perfectly account for the AGH approximation in (\ref
{sco1}) if the partial derivatives, at both the numerator and
denominator, of order greater than $2k$ ($\max m= k$) are zero. The
corresponding Laplace approximation can be rewritten by regrouping in
decreasing asymptotic order the elements that appear in both the
expansions, and by truncating the resulting series at an appropriate
point. In symbols,
%
%
\begin{equation}
\label{ratio} \frac{S(\bolds\theta; \hat{\mathbf
{z}})+c_{1}^{*}p^{-1}+ \cdots
+c_{r}^{*}p^{-r}+\cdots+c_{[k/3]}^{*}p^{-[k/3]}+ \RMO (p^{- [
{k}/{3}+1 ]} )}{1+c_{1}p^{-1}+ \cdots+c_{r}p^{-r}+\cdots
+c_{[k/3]}p^{-[k/3]}+ \RMO (p^{- [{k}/{3}+1 ]} )},
\end{equation}
where the coefficients $c_{r}, r=1,\ldots, [\frac{k}{3} ]$, are
given by
\[
c_{r}=\sum_{m=r+1}^{3r}
\frac{(-1)^{m-r}}{(2m)!}L_{p_{1}}(\hat{\mathbf{z}}) \cdots L_{p_{m-r}}(
\hat{\mathbf{z}})L^{q_{1}}(\hat{\mathbf{z}}) \cdots L^{q_{m}}(\hat{
\mathbf{z}})
\]
with $p_{1}|\cdots|p_{t}$ be a partition of $2m$ indices into $m-r$
blocks, each of size 3 or more, and $q_{1}|\cdots|q_{m}$ is a
partition of $2m$ indices into $m$ blocks, each of size 2. On the other
hand, the coefficients $c_{r}^{*}, r=1,\ldots, [\frac{k}{3} ]$,
results
\[
c_{r}^{*}=\sum_{m=r}^{3r}
\sum_{s=0}^{3r-m}\frac
{(-1)^{m-r}}{(2m)!}S_{j_{1},\ldots, j_{s}}(
\bolds\theta; \hat{\mathbf{z}})L_{p_{1}}(\hat{\mathbf{z}}) \cdots
L_{p_{m-r}}(\hat{\mathbf{z}})L^{q_{1}}(\hat{\mathbf{z}}) \cdots
L^{q_{m}}(\hat{\mathbf{z}}),
\]
where $p_{1}|\cdots|p_{m-r}$ is a partition of $2m-s$ indices into
$m-r$ blocks, each of size 3 or more, and $q_{1}|\cdots|q_{m}$ is a
partition of the same indices together with $\{j_{1},\ldots,j_{s}\}$
into $m$ blocks of size~2. 
Since $S_{j_{1},\ldots, j_{s}}(\bolds\theta; \hat{\mathbf
{z}})=\frac{\partial L_{j_{1},\ldots, j_{s}}(\hat{\mathbf
{z}})}{\partial\bolds\theta}$, all the first derivatives of the
score function will be zero due to the choice of $\hat{\mathbf{z}}$.

Based on long polynomial division, the approximated expected score
functions equivalent to (\ref{sco1}) are given by
%
%
\begin{equation}
\label{sco2}\sum_{l=1}^{n}
\tilde{E}_{\mathbf{z}|\mathbf
{y}} \bigl[S_{l}(\bolds\theta;
\mathbf{z}_{l}) \bigr]= \sum_{l=1}^{n}
\bigl[S_{l}(\bolds\theta; \hat{\mathbf{z}}_{l})+c_{1}^{**}p^{-1}+
\cdots+ c_{r}^{**}p^{-r}+ \cdots+ \RMO
\bigl(p^{- [{k}/{3}+1 ]} \bigr) \bigr],
\end{equation}
where the coefficients $c_{r}^{**}, r=1,\ldots, [\frac{k}{3} ]$,
can be determined as follows
\[
c_{r}^{**}=\bigl[c_{r}^{*}-S(
\bolds\theta; \hat{\mathbf{z}})c_{r}\bigr
]-c_{r-1}c_{1}^{*}-c_{r-2}c_{2}^{*}-
\cdots-c_{1}c_{r-1}^{*}
\]
being
\[
c_{r}^{*}-S(\bolds\theta; \hat{\mathbf{z}})c_{r}=
\sum_{m=r}^{3r}\sum
_{s=2}^{3r-m}\frac{(-1)^{m-r}}{(2m)!}S_{j_{1},\ldots,
j_{s}}(
\bolds\theta; \hat{\mathbf{z}})L_{p_{1}}(\hat{\mathbf{z}}) \cdots
L_{p_{m-r}}(\hat{\mathbf{z}})L^{q_{1}}(\hat{\mathbf{z}}) \cdots
L^{q_{m}}(\hat{\mathbf{z}})
\]
and, in particular, $c_{1}^{*}-S(\bolds\theta; \hat{\mathbf
{z}})c_{1}=\frac{1}{2}S_{j_{1},
j_{2}}L^{j_{1},j_{2}}(\hat{\mathbf{z}})$.

\subsection{Asymptotic behavior of the AGH-based estimators}\label{sec3.2}

To investigate the properties of the AGH approximated maximum
likelihood estimators $\hat{\bolds\theta}$, we analyze the
asymptotic behavior of the corresponding Laplace-based estimators
defined by (\ref{sco2}). Our arguments are similar to those of
Huber, Ronchetti and Victoria-Feser \cite{HubRonFe04}, who discussed classical
Laplace estimators in GLLVM, and Rizopoulos, Verbeke and Lesaffre \cite
{RizVer09}
who analyzed the consistency of fully exponential Laplace estimators in
joint models for survival and longitudinal data.
%
%
\begin{prop}[(Consistency)]\label{cons} Let $\bolds\theta
_{0} \in\Theta$ denote the true parameter value, then, under suitable
regularity conditions,
%
%
\begin{equation}
\label{consistency} (\hat{\bolds\theta} - \bolds\theta_{0})=\RMO_{p}
\bigl[\max\bigl(n^{-1/2},p^{- [{k}/{3}+1 ]} \bigr) \bigr].
\end{equation}
\end{prop}

Thus, $\hat{\bolds\theta}$ is consistent as long as
both $n$ and $p$ grow to $\infty$. A formal proof of Proposition \ref
{cons} is given in Appendix \ref{app2}. The $n^{-1/2}$ term comes from the
standard asymptotic theory, whereas the $p^{- [{k}/{3}+1
]}$ term derives from the AGH approximation. The requirement that $p$
grows to infinity is consistent with the fact that we are trying to
approximate the marginal density of each individual, that is, $f(\mathbf
{y}_{l}; \bolds\theta)$. 
However, in practical applications where $p$ and $k$ are both fixed,
the approximation error in the adaptive technique is $\RMO (p^{-
[{k}/{3}+1 ]} )$ as $n \rightarrow\infty$, and the
asymptotic properties of the AGH-based estimators should be evaluated
with respect to a perturbation of the true parameter $\bolds
\theta_{0}$.

For $k \geq3$, the AGH-based estimator is more accurate
than the classical $\RMO(p^{-1})$ Laplace-based estimators. Indeed, it
shares the same accuracy of higher (than one) order Laplace estimators,
but, with respect to this latter, the adaptive Gauss--Hermite is easier
to be implemented, since it avoids derivative computations.

Based on the derivation of (\ref{consistency}) as presented in the
Appendix \ref{app2}, we can deduce that, if $p=\RMO(n^{\rho})$ for $\rho>\frac{1}{
[{k}/{3}+1 ]}$, then the AGH-based estimators will be
asymptotically equivalent to the true maximum likelihood estimators
that solve $S(\bolds\theta)=0$. However, in general, they are not
maximum likelihood estimators because of the approximation, but, as
discussed by Huber, Ronchetti and Victoria-Feser \cite{HubRonFe04} for
classical Laplace estimators in GLLVM, they belong to the class of
$M$-estimators. The latter are implicitly defined through a general
$\Psi$-function as the solution in $\bolds\theta$ of
\[
\sum_{l=1}^{n}\Psi(\mathbf{y}_{l};
\bolds\theta)=0.
\]
The $\Psi$-function for the AGH-based estimators are given by (\ref{sco1}).

%
\begin{prop}[(Asymptotic normality)] If $\bolds\theta
_{0}$ is an interior point of the parameter space $\Theta$ and
$B(\bolds\theta_{0})=-E [\frac{\partial\Psi(\mathbf
{y}_{l};\bolds\theta_{0})}{\partial\bolds\theta}
]=-E [\frac{\partial^{2} \tilde{\ell}(\bolds\theta
_{0})}{\partial\bolds\theta\,\partial\bolds\theta^{T}}
]$ is nonsingular, the AGH-based estimators are asymptotically normal,
that is,
%
%
\begin{equation}
\label{Mest}\sqrt{n}(\hat{\bolds\theta} - \bolds\theta_{0})
\rightarrow^{D} \mathit{MVN}\bigl(\mathbf{0}, B(\bolds
\theta_{0})^{-1}A(\bolds\theta_{0})\bigl[B(
\bolds\theta_{0})^{-1}\bigr]^{T}\bigr)
\end{equation}
with $A(\bolds\theta_{0})=E[\Psi(\mathbf{y}_{l};\bolds\theta
_{0})\Psi^{T}(\mathbf{y}_{l};\bolds\theta_{0})]=E [\frac
{\partial\tilde{\ell}(\bolds\theta_{0})}{\partial\bolds
\theta}\,{\frac{\partial\tilde{\ell}(\bolds\theta_{0})}{\partial
\bolds\theta}}^{T} ]$.
\end{prop}

The regularity conditions that ensure consistency and
asymptotic normality of the AGH-based $M$-estimators have to be checked
for the particular conditional distribution of each $y_{j}$ (Huber,
Ronchetti and Victoria-Feser \cite{HubRonFe04}). For classical Laplace-based
estimators, Huber, Scaillet and Victoria-Feser \cite{HubScaFe09}
analyzed these
conditions for ordered multinomial distributed manifest variables. A
formal derivation for the $M$-estimators discussed here in the case of
binary observed variables is provided in Appendix \ref{app3}.

\section{Monte Carlo simulations}\label{sec4}
In this section, we investigate empirically the finite sample
performance of the adaptive Gauss--Hermite and related Laplace-based
estimators. We focus on latent variable models for binary data, since
in this case the differences between numerical techniques should be
better highlighted
(Joe \cite{Joe08}). We consider two simulation scenarios
characterized by an increasing number of observed and latent variables.
In particular, we generate data from a population that consists of six
items satisfying a three factor model, and from a population based on
ten observed variables that satisfy a five factor model. In both cases,
the population parameters have been chosen in such a way that the
item-specific intercepts and the factor loadings are drawn randomly
from a log-normal distribution, with some loadings fixed to 0 to get
unique solutions. For each scenario, 100 random samples have been
considered with 200 subjects.

A crucial choice in the application of the AGH quadrature is the number
of points needed to adequately approximate the likelihood function. In
the simulation study, we follow Schilling and Bock \cite{SchBo05} who
suggested to select, in presence of binary data, five and three
quadrature points for the three and five factor model,\vadjust{\goodbreak} respectively. In
both cases, the performance of AGH is compared with that of the Laplace
approximation of order $\RMO(p^{-2})$. The estimation is performed through
the direct maximization algorithm described in Section~\ref{sec2}, whose
mathematical details for the case of binary observed items are provided
in Appendix \ref{app3}. The algorithm is written in the statistical language R
(R Development Core Team \cite{R10}) and the program is available from the
authors on request.

In the case of each simulation, the true values used to generate the
samples, the mean values of the estimated parameters across
simulations, together with their corresponding standard deviations
obtained from the simulated results, the mean estimated standard errors
obtained from (\ref{Mest}) and the Root Mean Square Error (RMSE) are
reported. Furthermore, in order to better highlight the computational
burden of each technique under the different conditions of study, we
report the average (over all the generated samples) computational time
in minutes (Avg min) and the average number of iterations (Avg iter)
required by the algorithm to get the convergence in a sample (obtain on
Intel Core i7 quad-core, 3.1 GHz CPU with 16 Gb RAM).

%
%
\begin{table}
\caption{True values, mean, simulated standard deviations (S.D.), root
mean square error (RMSE) and estimated standard errors (S.E.) of the
parameter estimates for AGH based on 5 quadrature points and for second
order Laplace ($Lap2$) approximation
in data generated by a three factor model with six ($p=6$) items
observed on $n=200$ subjects}
\label{Tab1}
\begin{tabular*}{\tablewidth}{@{\extracolsep{\fill}}llllllllll}
\hline
& \multicolumn{4}{l}{\emph{AGH}} && \multicolumn{4}{l@{}}{\emph{Lap2}}
\\[-4pt]
& \multicolumn{4}{l}{\hrulefill} && \multicolumn{4}{l@{}}{\hrulefill}
\\
True & Mean & S.D. & RMSE&S.E. && Mean & S.D.& RMSE&S.E. \\
\hline
$\alpha_{11}= 1.01$ & 0.72 & 0.37 & 0.47&0.19 && 1.32& 0.70 & 0.77&0.19
\\
$\alpha_{21}= 0.91$ & 1.17 & 0.36 & 0.45&0.46 && 1.11& 0.38 &0.42&0.49
\\
$\alpha_{31}= 0.50$ & 0.39 & 0.31 & 0.34&0.16 && 0.35&0.34 & 0.37&0.36
\\
$\alpha_{41}= 0.74$ & 0.99 & 0.38 & 0.45&0.27 && 1.14&0.18 & 0.44
&0.25\\
$\alpha_{51}= 1.16$ & 1.39 & 0.37 & 0.44&0.57 &&1.68& 0.37 & 0.64&0.66
\\
$\alpha_{61}= 1.22$ & 1.54 & 0.44 & 0.55&0.26 &&1.23& 0.52 & 0.52&0.42
\\
$\alpha_{12}= 0.00$ & -- & -- & -- &-- &&-- & -- & -- &--\\
$\alpha_{22}= 0.83$ & 0.45 & 0.30 & 0.49&0.32 && 0.21&0.69&0.93&0.54 \\
$\alpha_{32}= 0.44$ & 1.02 & 0.38& 0.69&0.42 && 1.06&0.36 &0.71&0.61 \\
$\alpha_{42}= 0.88$ & 1.15 & 0.42 & 0.50&0.57 && 1.13& 0.48 &0.54&0.53
\\
$\alpha_{52}= 1.73$ & 2.54 & 0.53 & 0.96&0.91 && 2.53& 0.37 & 0.88&0.85
\\
$\alpha_{62}= 1.46$ & 1.43 & 0.52 & 0.52&0.46 && 1.36& 0.73 & 0.73&0.53
\\
$\alpha_{13}= 0.00$ & -- & -- & -- &-- &&-- & -- & -- &-- \\
$\alpha_{23}= 0.00$ & -- & -- & --&-- && --&--&-- &-- \\
$\alpha_{33}= 1.45$ & 1.08 & 0.44& 0.58&0.49 && 1.21&0.68 &0.72&0.69 \\
$\alpha_{43}= 1.05$ & 1.52 & 0.42 & 0.64&0.62 && 1.49&0.36 &0.57&0.47 \\
$\alpha_{53}= 0.62$ & 0.93 & 0.36 & 0.48&0.58 && 0.80& 0.37 & 0.41&0.66
\\
$\alpha_{63}= 0.91$ & 0.98 & 0.42 & 0.43 &0.42 && 0.53&0.34 & 0.51&0.41
\\[6pt]
Avg iter&\multicolumn{4}{l}{9.21}&&\multicolumn{4}{l@{}}{324.94}\\
Avg min&\multicolumn{4}{l}{$3'\,52''$}&&\multicolumn{4}{l@{}}{\hphantom{0}$24'\,26''$}\\
\hline
\end{tabular*}
\end{table}

Table~\ref{Tab1} studies the performance of the AGH based on five
quadrature points and of the second order Laplace approximation on the
data generated by the three factor model. The results show that the two
techniques provide similar RMSE values for almost all the model
parameters. Indeed, even if the Laplace seems to introduce a slightly
larger bias in some estimates than the AGH, the simulated standard
deviations, that are a measure of the sampling variability of the
estimated parameters, are quite close. The estimated standard errors
are generally larger than the simulated standard deviations for both
the techniques, and closer to the corresponding RMSE. As expected, the
main difference between the two techniques is computational. The
algorithm based on the adaptive quadrature achieves convergence for a
sample, on average, in ten iterations, that is in less than four
minutes, whereas the second order Laplace requires, on average, more
than 300 iterations to get the convergence in a sample, that means
almost thirty minutes. This is more evident in Table~\ref{Tab2} that
%
%
\begin{table}
\caption{True values, mean, simulated standard deviations (S.D.), root
mean square error (RMSE) and estimated standard errors (S.E.) of the
parameter estimates for AGH based on three quadrature points, and for
second order Laplace ($Lap2$) approximation
in the data generated by a five factor model with ten ($p=10$) items
observed on $n=200$ individuals}
\label{Tab2}
\begin{tabular*}{\tablewidth}{@{\extracolsep{\fill}}llllllllll@{}}
\hline
& \multicolumn{4}{l}{\emph{AGH}} && \multicolumn{4}{l@{}}{\emph{Lap2}} \\[-4pt]
& \multicolumn{4}{l}{\hrulefill} && \multicolumn{4}{l@{}}{\hrulefill}
\\
True & Mean & S.D. & RMSE &S.E. && Mean & S.D. & RMSE &S.E. \\
\hline
$\alpha_{11}= 1.01$ & 0.70 & 0.34 & 0.46&0.64 && 1.28 & 0.46 &0.53&
0.54 \\
$\alpha_{21}= 0.91$ & 1.27 & 0.27& 0.45&0.56 && 1.33 & 0.51& 0.66 &0.56
\\
$\alpha_{31}= 0.50$ & 0.87 & 0.41 & 0.55 &0.52 && 0.57 & 0.22 & 0.23
&0.42 \\
$\alpha_{41}= 0.74$ & 1.14 & 0.39 & 0.56&0.61 && 0.85 & 0.33 & 0.35
&0.58 \\
$\alpha_{51}= 1.16$ & 1.83 & 0.26 & 0.72 &0.71 && 1.98 & 0.48 & 0.95
&0.81 \\
$\alpha_{61}= 1.22$ & 1.22 & 0.39 & 0.39&0.42 && 0.66 & 0.23 &0.60&0.62
\\
$\alpha_{71}= 0.55$ & 0.48 & 0.31& 0.32&0.27 && 0.59 & 0.22& 0.22&0.27
\\
$\alpha_{81}= 0.83$ & 1.10 & 0.35 & 0.44&0.45 && 0.90 & 0.26 & 0.27
&0.35 \\
$\alpha_{91}= 0.44$ & 1.01 & 0.28 & 0.63 &0.64 && 1.05 & 0.21 & 0.65
&0.64 \\
$\alpha_{101}= 0.88$ & 1.05 & 0.30 & 0.35 &0.36 && 1.17 & 0.18 & 0.34
&0.43 \\
$\alpha_{12}= 0.00$ & -- & -- & -- &-- && -- & -- & -- &-- \\
$\alpha_{22}= 1.46$ & 1.39 & 0.26 & 0.26&0.26 && 1.21 & 0.44 & 0.50 &
0.28 \\
$\alpha_{32}= 0.89$ & 0.59 & 0.50& 0.58&0.57 && 0.80 & 0.38 & 0.39
&0.53 \\
$\alpha_{42}= 1.64$ & 1.27 &0.23 & 0.44 &0.45 && 1.37 & 0.38 & 0.47
&0.52 \\
$\alpha_{52}= 1.45$ & 0.60 &0.35 & 0.92 &0.91 && 0.59 & 0.46 &0.98
&0.91 \\
$\alpha_{62}= 1.05$ & 0.92 & 0.37 & 0.39 &0.38 && 1.06 & 0.19 &0.19
&0.28 \\
$\alpha_{72}= 0.62$ & 0.68 & 0.42 & 0.42&0.43 && 0.80 & 0.34 &0.39&0.43
\\
$\alpha_{82}= 0.91$ & 0.40 &0.32& 0.60 &0.58 && 0.19 & 0.41 &0.83 &0.60
\\
$\alpha_{92}= 1.59$ & 1.22 &0.31 & 0.48 &0.48 && 2.02 & 0.59 &0.73
&0.52 \\
$\alpha_{102}= 1.27$ & 0.95 &0.32 & 0.46&0.46 && 1.22 & 0.27 &0.27&0.39
\\
$\alpha_{13}= 0.00$ & -- & -- & -- &-- && -- & -- & -- &-- \\
$\alpha_{23}= 0.00$ & -- & -- & --&-- && -- & -- & -- &-- \\
$\alpha_{33}= 0.71$ & 1.10 & 0.45& 0.59&0.62 && 1.13 & 0.51 & 0.66 &
0.65 \\
$\alpha_{43}= 0.35$ & 1.02 &0.29 & 0.73 &0.74 && 1.01 & 0.15 & 0.68
&0.65 \\
$\alpha_{53}= 0.53$ & 1.46 &0.28 & 0.97 &0.98 && 1.46 & 0.16 &
0.95&0.98 \\
\hline
\end{tabular*}
\end{table}
shows the results for the five factor model. In this specific case, the
adaptive Gauss--Hermite has been applied with three quadrature points.
For this latter, the algorithm requires, on average, less than ten
iterations to get convergence in a sample, that is less than five
minutes, whereas the algorithm based on the second order Laplace
approximation is much slower than in the case of the three factor
model. As before, it reaches convergence, on average, in almost 350
iterations, but now it requires almost four hours to obtain the
solution for one sample. However, the estimates derived by applying the
two techniques are quite comparable in terms of bias, standard
deviations and RMSE, with similar conclusions to those drawn for the
first scenario.\looseness=-1

\setcounter{table}{1}
\begin{table}
\caption{(Continued)}
\begin{tabular*}{\tablewidth}{@{\extracolsep{\fill}}llllllllll@{}}
\hline
& \multicolumn{4}{l}{\emph{AGH}} && \multicolumn{4}{l@{}}{\emph{Lap2}} \\[-4pt]
& \multicolumn{4}{l}{\hrulefill} && \multicolumn{4}{l@{}}{\hrulefill}
\\
True & Mean & S.D. & RMSE &S.E. && Mean & S.D. & RMSE &S.E. \\
\hline
$\alpha_{63}= 0.83$ & 0.97 & 0.40& 0.42&0.45 && 0.64 & 0.25 &0.31&0.41
\\
$\alpha_{73}= 0.71$ & 1.12 & 0.36 & 0.55 &0.56&& 1.52 & 0.37 & 0.69
&0.59 \\
$\alpha_{83}= 0.65$ & 1.36 & 0.30& 0.77&0.77 && 1.27 & 0.47 &0.78 &0.77
\\
$\alpha_{93}= 0.95$ & 1.19 &0.30& 0.39&0.41 && 0.84 & 0.17 & 0.21 &0.39
\\
$\alpha_{103}= 0.88$ & 1.23 &0.38 & 0.52&0.54 && 0.68 & 0.18 &0.37&0.39
\\
$\alpha_{14}= 0.00$ & -- & -- & -- &-- && -- & -- & -- &-- \\
$\alpha_{24}= 0.00$ & -- & -- & --&-- && -- & -- & -- &-- \\
$\alpha_{34}= 0.00$ & -- & --& --&-- && -- & -- & -- &--\\
$\alpha_{44}= 1.10$ & 1.42 &0.37& 0.49 &0.51 && 1.95 & 0.40 & 0.65&0.54\\
$\alpha_{54}= 0.50$ & 0.84 &0.29 & 0.45&0.46 && 0.95 & 0.57 &0.53 &0.49
\\
$\alpha_{64}= 0.49$ & 0.93 & 0.42 & 0.61 &0.62 && 0.05 & 0.61 &0.56
&0.62 \\
$\alpha_{74}= 1.20$ & 0.67 & 0.51 & 0.73&0.74 && 0.50 & 0.24 & 0.74
&0.74 \\
$\alpha_{84}= 0.41$ & 0.43 & 0.36& 0.36 &0.36 && 0.11 & 0.31 &0.43&0.38
\\
$\alpha_{94}= 0.85$ & 0.82 &0.37 & 0.37 &0.41&& 0.43 & 0.20 &0.47&0.40
\\
$\alpha_{104}= 0.72$ & 1.03 &0.37 & 0.48&0.48 && 1.11 & 0.21&0.44&0.42
\\
$\alpha_{15}= 0.00$ & -- & -- & -- &-- && -- & -- & -- &-- \\
$\alpha_{25}= 0.00$ & -- & -- & -- &--&& -- & -- & -- &-- \\
$\alpha_{35}= 0.00$ & -- & --& --&-- && -- & -- & -- &--\\
$\alpha_{45}= 0.00$ & -- &-- & -- &-- && -- & -- & --&-- \\
$\alpha_{55}= 0.62$ & 0.80 &0.28& 0.33&0.34 && 0.53 & 0.35 & 0.36& 0.34
\\
$\alpha_{65}= 0.99$ & 1.32 & 0.32 & 0.47 &0.48 && 1.23 & 0.62 & 0.57
&0.52 \\
$\alpha_{75}= 1.12$ & 1.04 & 0.36 & 0.37 &0.38&& 0.97 & 0.36 &0.35&0.40
\\
$\alpha_{85}= 0.86$ & 1.05 & 0.40& 0.45&0.46 && 1.42 & 0.52 &0.37&0.43
\\
$\alpha_{95}= 0.71$ & 0.67 &0.33 & 0.33 &0.36&& 0.94 & 0.22 &0.32&0.35
\\
$\alpha_{105}= 1.39$ & 1.32 &0.35 & 0.36 &0.35 && 1.73 & 0.54
&0.25&0.35 \\[6pt]
Avg iter&\multicolumn{4}{l}{9.15}&&\multicolumn{4}{l@{}}{344.27}\\
Avg min&\multicolumn{4}{l}{$4'\,47''$}&&\multicolumn{4}{l@{}}{$239'\,16''$}\\
\hline
\end{tabular*}
\end{table}

To better investigate the properties of the adaptive ML estimators, a
further simulation study has been conducted in order to understand how
much contribution is due to the approximation error and how much is due
to the term $\RMO(n^{1/2})$ in the rate of consistency (\ref
{consistency}). For the three factor model, the performance of the
adaptive quadrature based on five quadrature points ($\mathit{AGH}5$) has been
analyzed in presence of small ($n=200$) and large ($n=1000$) samples.
Furthermore, for the smallest sample size, the behavior of AGH has been
studied by also considering nine ($\mathit{AGH}9$) and fifteen ($\mathit{AGH}15$)
quadrature points. All the results are illustrated in Table~\ref{Tab3}.
We recall that $\mathit{AGH}9$ shares the same asymptotic properties of the
Laplace estimator of order $\RMO(p^{-4})$, whereas when fifteen quadrature
points are used, the error rate is of order $\RMO(p^{-6})$. However, the
performance of these Laplace estimators is not analyzed since they
require a lot of time just to run a simple simulation example as the
one considered here.

%
\begin{table}
\caption{True values, mean, simulated standard deviations (S.D.), root
mean square error (RMSE) and estimated standard errors (S.E.) of the
parameter estimates for AGH based on 5 ($\mathit{AGH}5$), 9 ($\mathit{AGH}9$) and 15
($\mathit{AGH}15$) quadrature points in data generated by a three factor model
with six observed items ($p=6$)}
\label{Tab3}
\begin{tabular*}{\tablewidth}{@{\extracolsep{\fill}}llllllllll@{}}
\hline
& \multicolumn{9}{l@{}}{\emph{AGH5}} \\[-4pt]
& \multicolumn{9}{l@{}}{\hrulefill}\\
& \multicolumn{4}{l}{$n=200$} && \multicolumn{4}{l@{}}{$n=1000$} \\[-4pt]
& \multicolumn{4}{l}{\hrulefill} && \multicolumn{4}{l@{}}{\hrulefill} \\
True & Mean & S.D. & RMSE&S.E. && Mean & S.D.& RMSE&S.E. \\
\hline
$\alpha_{11}= 1.01$ & 0.72 & 0.37 & 0.47&0.19 && 0.89&0.12&0.17&0.17 \\
$\alpha_{21}= 0.91$ & 1.17 & 0.36 & 0.45&0.46 && 0.98& 0.06 &0.10&0.19
\\
$\alpha_{31}= 0.50$ & 0.39 & 0.31 & 0.34&0.16 && 0.53&0.00 & 0.03&0.08
\\
$\alpha_{41}= 0.74$ & 0.99 & 0.38 & 0.45&0.27 && 0.90&0.00 & 0.16
&0.15\\
$\alpha_{51}= 1.16$ & 1.39 & 0.37 & 0.44&0.57 &&1.35& 0.00 & 0.19&0.19
\\
$\alpha_{61}= 1.22$ & 1.54 & 0.44 & 0.55&0.26 &&1.31& 0.22 & 0.24&0.26
\\
$\alpha_{12}= 0.00$ & -- & -- & -- &-- &&-- & -- & -- &--\\
$\alpha_{22}= 0.83$ & 0.45 &0.30 & 0.49&0.32 && 0.73&0.00&0.11&0.12 \\
$\alpha_{32}= 0.44$ & 1.02 & 0.38& 0.69&0.42 && 0.87&0.00&0.42&0.31 \\
$\alpha_{42}= 0.88$ & 1.15 & 0.42 & 0.50&0.57 && 1.07& 0.00 &0.19&0.17
\\
$\alpha_{52}= 1.73$ & 2.54 &0.53 & 0.96&0.91 && 2.10& 0.62 & 0.75&0.71
\\
$\alpha_{62}= 1.46$ & 1.43 &0.52 & 0.52&0.46 && 1.45& 0.28 & 0.28&0.26
\\
$\alpha_{13}= 0.00$ & -- & -- & -- &-- &&-- & -- & -- &-- \\
$\alpha_{23}= 0.00$ & -- & -- & --&-- && --&--&-- &-- \\
$\alpha_{33}= 1.45$ & 1.08 & 0.44& 0.58&0.49 && 1.19&0.31 &0.41&0.39 \\
$\alpha_{43}= 1.05$ & 1.52 & 0.42 & 0.64&0.62 && 1.27&0.38 &0.45&0.44 \\
$\alpha_{53}= 0.62$ & 0.93 &0.36 & 0.48&0.58 && 0.80&0.35 & 0.39&0.38 \\
$\alpha_{63}= 0.91$ & 0.98 &0.42 & 0.43&0.42 && 0.90&0.29 & 0.35&0.32\\
[6pt]
Avg iter&\multicolumn{4}{l}{9.21}&&\multicolumn{4}{l@{}}{11.66}\\
Avg min&\multicolumn{4}{l}{$3'\,52''$}&&\multicolumn{4}{l@{}}{$41'\,28''$}\\
\hline
\end{tabular*}    \vspace*{-3pt}
\end{table}

\setcounter{table}{2}
\begin{table}
\caption{(Continued)}
\begin{tabular*}{\tablewidth}{@{\extracolsep{\fill}}llllllllll@{}}
\hline
& \multicolumn{4}{l}{\emph{AGH9}} && \multicolumn{4}{l@{}}{\emph{AGH15}}\\[-4pt]
& \multicolumn{4}{l}{\hrulefill} && \multicolumn{4}{l@{}}{\hrulefill}\\
& \multicolumn{4}{l}{$n=200$}&& \multicolumn{4}{l@{}}{$n=200$}\\[-4pt]
& \multicolumn{4}{l}{\hrulefill} && \multicolumn{4}{l@{}}{\hrulefill}\\
True & Mean & S.D. & RMSE&S.E. && Mean & S.D.& RMSE&S.E. \\
\hline
$\alpha_{11}= 1.01$ & 0.82&0.32&0.37&0.17&&0.83&0.28&0.31&0.29 \\
$\alpha_{21}= 0.91$ &1.01&0.33&0.34&0.44&&1.00&0.24&0.26&0.35\\
$\alpha_{31}= 0.50$ &0.39&0.32&0.34&0.38&&0.60&0.24&0.25&0.28 \\
$\alpha_{41}= 0.74$ &0.94&0.40&0.45&0.28&&0.93&0.32&0.32&0.27\\
$\alpha_{51}= 1.16$ &1.29&0.36&0.38&0.54&&1.38&0.32&0.39&0.39\\
$\alpha_{61}= 1.22$ &1.56&0.35&0.49&0.30&&1.48&0.41&0.44&0.32 \\
$\alpha_{12}= 0.00$ &--&--&--&--&&--&--&--&--\\
$\alpha_{22}= 0.83$ & 0.55&0.31&0.42&0.43&&0.63&0.24&0.27&0.32\\
$\alpha_{32}= 0.44$ &0.93&0.44&0.65&0.51&&0.89&0.18&0.42&0.43 \\
$\alpha_{42}= 0.88$ &1.04&0.42&0.45&0.41&&1.06&0.25&0.26&0.31 \\
$\alpha_{52}= 1.73$ &2.40&0.45&0.80&0.77 &&2.22&0.50&0.70&0.70\\
$\alpha_{62}= 1.46$ &1.43&0.50&0.50&0.33&&1.43&0.48&0.48&0.43\\
$\alpha_{13}= 0.00$ &--&--&--&--&&--&--&--&--\\
$\alpha_{23}= 0.00$ &--&--&--&--&&--&--&--&-- \\
$\alpha_{33}= 1.45$ &1.10&0.36&0.50&0.51&&1.12&0.40&0.50&0.47\\
$\alpha_{43}= 1.05$ &1.46&0.37&0.56&0.52&&1.42&0.33&0.49&0.46\\
$\alpha_{53}= 0.62$ &0.94&0.40&0.51&0.49&&0.89&0.41&0.45&0.45\\
$\alpha_{63}= 0.91$ &0.99&0.43&0.43&0.42&&0.95&0.45&0.46&0.46\\[6pt]
Avg iter&\multicolumn{4}{l}{\hphantom{0}9.40}&&\multicolumn{4}{l@{}}{\hphantom{0}12.74}\\
Avg min&\multicolumn{4}{l}{$19'\,37''$}&&\multicolumn{4}{l@{}}{$107'\,11''$}\\
\hline
\end{tabular*}
\end{table}

In Table~\ref{Tab3}, it can be noticed that $\mathit{AGH}5$ performs better in
the largest sample than in the smallest one in terms of both bias and RMSE. The estimated
standard errors for all the parameters become smaller and also closer
to the root mean square error as the sample size increases.

On the other hand, increasing the number of quadrature
points, the AGH performs better in terms of bias and RMSE with slight
differences between $\mathit{AGH}9$ and $\mathit{AGH}15$, mainly due to a less
variability in the estimates for the latter than for the former.
However, the adaptive Gauss--Hermite quadrature becomes more
computational intensive as the number of quadrature points increases.
As shown in Table~\ref{Tab3}, the algorithm needs, on average, almost
ten iterations to get convergence in presence of both five and nine
quadrature points. However, whereas in the former case, the solution
for a sample is obtained, on average, in four minutes, more than
fifteen minutes are required in the latter case. This is more evident
for $\mathit{AGH}15$, for which the algorithm gets convergence, on average, in
almost thirteen iterations, but requiring more than one hour and a half
to obtain the solution for one sample. It is also evident that to get
the same accuracy in the estimates observed for $\mathit{AGH}5$ in the largest
sample, in presence of binary data, more than fifteen quadrature points
per dimension should be considered in small samples.\vspace*{-3pt}

\section{Conclusions}\label{sec5}\vspace*{-3pt}

In this paper, we have investigated the theoretical properties of
adaptive Gauss--Hermite based estimators in the GLLVM framework.
Recently, the adaptive quadrature has played a prominent role in the
latent variable model literature for approximating integrals defined
over the latent space. It allows to overcome the main limitations of
the commonly used techniques, such as the Gauss--Hermite quadrature and
the standard Laplace approximation. Indeed, AGH is applicable to
problems\vadjust{\goodbreak} involving high-dimensional integrals where the former becomes
impractical or computationally intensive, and it provides more accurate
estimates than the latter, particularly when used for binary or ordinal
data with small sample sizes
(Joe \cite{Joe08}).

We have proved that, for multidimensional integrals, the AGH solution
is asymptotically equivalent to the Laplace approximation that involves
specific higher (than two) order derivatives of the integrand. Higher
order Laplace approximations have been suggested in several papers on
generalized linear models
(Raudenbush, Yang and Yosef \cite{Rade00}, Evangelou, Zhu and Smith
\cite{Eva11},
Bianconcini and Cagnone \cite{BiaCag12}) as an
alternative to classical methods for improving the accuracy of the
estimates. This extension has been motivated by the well-known
asymptotic properties that characterize the Laplace method, and by the
fact that the approach does not suffer from the ``curse of
dimensionality''. However, the inclusion of higher order terms is
computationally demanding as the order of the approximation increases.
On the other hand, the AGH quadrature is easier to be implemented, but
of course its computational complexity increases as the number of
latent variables increases. Hence, AGH and higher order Laplace
approximations can be seen as complementary approaches that share the
same asymptotic properties.

We have shown that the AGH-based estimators are consistent as the
sample size and number of observed variables grow to infinity. The
convergence rate of these estimators depends also on the number of
quadrature points used for each dimension. In general, these estimators
are less efficient than maximum likelihood estimators because of the
approximation, but belong to the class of $M$-estimators, for which the
asymptotic properties are well-known such that correct inference can be
performed.

%
\begin{appendix}\label{app}
\section{Asymptotic behavior of the multivariate AGH~approximation}\label{app1}
The higher order Laplace approximation of (\ref{eq11}) is derived
by considering
\[
f(\mathbf{y}; \bolds\theta)=\int_{\mathbb{R}^{q}}\RMe^{-L(\mathbf
{z})}\mrmd
\mathbf{z},
\]
where $L(\mathbf{z})=-[\log g(\mathbf{y}|\mathbf{z}; \bolds\theta
) + \log h(\mathbf{z})]$, being $L(\mathbf{z})=\RMO(p)$. It is based on
the Taylor series expansion of $L$ around its minimum $\hat{\mathbf
{z}}$, that is,
%
%
\begin{equation}
\label{logtay} L(\mathbf{z})=L(\hat{\mathbf{z}}) + \frac{1}{2}L_{i_{1},i_{2}}(
\hat{\mathbf{z}}) (\mathbf{z}-\hat{\mathbf{z}})^{i_{1},i_{2}}+\sum
_{m=3}^{\infty}\frac{1}{m!}L_{i_{1},\ldots,i_{m}}(\hat{
\mathbf{z}}) (\mathbf{z}-\hat{\mathbf{z}})^{i_{1},\ldots,i_{m}}.
\end{equation}
Substituting (\ref{logtay}) into the integral, we obtain
\begin{eqnarray*}
&&
\int_{\mathbb{R}^{q}}\exp\bigl[-L(\mathbf{z})\bigr]\mrmd
\mathbf{z}\\[-1pt]
&&\quad=(2\pi)^{q/2}|\bolds\Psi|^{{1/2}}\RMe^{-L(\hat{\mathbf{z}})}
\int_{\mathbb{R}^{q}}h_{1}(\mathbf{z};\hat{\mathbf{z}},
\bolds\Psi)\exp\Biggl[\sum_{m=3}^{\infty}
\frac{(-1)}{m!}L_{i_{1},\ldots,i_{m}}(\hat{\mathbf{z}}) (\mathbf
{z}-\hat{
\mathbf{z}})^{i_{1},\ldots,i_{m}} \Biggr]\mrmd \mathbf{z}
\\[-1pt]
&&\quad=(2\pi)^{q/2}|\bolds\Psi|^{{1/2}}\RMe^{-L(\hat
{\mathbf{z}})}
\Biggl[1-\sum_{m=2}^{\infty}\sum
_{P,Q}\frac
{(-1)^{t}}{(2m)!}L_{p_{1}}(\hat{\mathbf{z}})
\cdots L_{p_{t}}(\hat{\mathbf{z}})L^{q_{1}}(\hat{\mathbf{z}})
\cdots L^{q_{m}}(\hat{\mathbf{z}}) \Biggr],
\end{eqnarray*}
where the second sum is over all partitions $P,Q$, such that
$P=p_{1}|\cdots|p_{t}$ is a partition of $2m$ indices into $t$ blocks,
each of size 3 or more, and $Q=q_{1}|\cdots|q_{m}$ is a partition of
$2m$ indices into $m$ blocks, each of size 2. Each component
$L^{q_{k}}(\hat{\mathbf{z}}), k=1,\ldots, m$, refers to specific
elements of the covariance matrix $\bolds\Psi$.

We want here to show that the exact higher order Laplace solution for
the integral (\ref{eq11}) is equivalent to the one based on the AGH
quadrature given in (\ref{taylor2}). To do so, we need to show that
%
%
\begin{eqnarray}
\label{equivalence}
&&1+\sum_{m=3}^{\infty}\sum
_{P}(-1)^{t}L_{p_{1}}(\hat{
\mathbf{z}}) \cdots L_{p_{t}}(\hat{\mathbf{z}}) (\mathbf{z}- \hat{
\mathbf{z}})^{i_{1},\ldots, i_{m}}\nonumber\\[-8pt]\\[-8pt]
&&\quad= \sum_{m=3}^{\infty}
c_{i_{1},\ldots, i_{m}}(\hat{\mathbf{z}}) (\mathbf{z}- \hat{\mathbf
{z}})^{i_{1},\ldots, i_{m}}.\nonumber
\end{eqnarray}
At this regard, we can notice that, based on the exlog relations, the
LHS term of (\ref{equivalence}) is equal to $\exp[\sum
_{m=3}^{\infty} L_{i_{1},\ldots, i_{m}}(\hat{\mathbf{z}})(\mathbf{z}-
\hat{\mathbf{z}})^{i_{1},\ldots, i_{m}} ]$. This higher order
term can be rewritten as
\[
\pi(\mathbf{z})=\exp\bigl[-L(\mathbf{z})+L(\hat{\mathbf{z}})+
\tfrac
{1}{2}L_{i_{1},i_{2}}(\hat{\mathbf{z}}) (\mathbf{z}-\hat{\mathbf
{z}})^{i_{1},i_{2}} \bigr]
\]
that it is equal to
\[
\frac{(2\pi)^{-{q}/{2}}|\bolds\Psi|^{-{1}/{2}}g(\mathbf
{y}|\mathbf{z}; \bolds\theta)h(\mathbf{z})}{g(\mathbf{y}|\hat
{\mathbf{z}};\bolds\theta)h(\hat{\mathbf{z}})h_{1}(\mathbf{z};\hat
{\mathbf{z}}, \bolds\Psi)}=\frac{\nu(\mathbf{z})}{\nu(\hat{\mathbf
{z}})}=c(\mathbf{z}).
\]
Hence, the Taylor series expansion of $\pi(\mathbf{z})$ around the
minimum $\hat{\mathbf{z}}$ can be written as
\begin{eqnarray*}
\pi(\mathbf{z})&=&1+\sum_{m=3}^{\infty}
\frac{1}{m!}\pi_{i_{1},\ldots,i_{m}}(\hat{\mathbf{z}}) (\mathbf{z}-\hat{
\mathbf{z}})^{i_{1},\ldots,i_{m}}
\\
&=&1+\sum_{m=3}^{\infty}
\frac{1}{m!}c_{i_{1},\ldots,i_{m}}(\hat{\mathbf{z}}) (\mathbf{z}-\hat{
\mathbf{z}})^{i_{1},\ldots,i_{m}}.
\end{eqnarray*}
It follows that the AGH solution (\ref{taylor2}) and the Laplace one
(\ref{Ei2}) are equivalent. Based on this relationship, it is possible
to derive the asymptotic error associated with the AGH approximation
(\ref{taylor3}) evaluating the equivalent Laplace approximation
obtained by truncating (\ref{Ei2}) at $m=k$. Shun and McCullagh \cite
{ShuMc95} proved that,
for fixed $q$, the usual asymptotic order of the term corresponding to
the bipartition $(P,Q)$ is $\RMO(p^{t-m})$. The error rate of the AGH
based on $k$ quadrature points is the same associated to the
bipartition $(P,Q)$ of $2(k+1)$ indices in the expansion (\ref{Ei2}).
In this case, the maximum number of blocks, each of size at least 3,
for $2(k+1)$ indices is $ [\frac{2(k+1)}{3} ]$, where $[r]$
indicates the largest integer not exceeding $r$. Hence, being $m=k+1$,
the AGH based on $k$ quadrature points has associated asymptotic order
equal to $\RMO (p^{- [{k}/{3}+1 ]} )$.

\section{Consistency of the AGH-based estimators}\label{app2}

This section concerns with the consistency of the AGH-based estimators.
All the following proofs proceed along the lines of Vonesh
\cite{Von96}, who derived the rate of convergence of the estimator
based on the classical Laplace approximation for nonlinear mixed effect
models, and of Rizopoulos, Verbeke and Lesaffre \cite{RizVer09}, who
derived that
rate for fully exponential Laplace based estimators in joint models for
longitudinal and survival data. In particular, we work under the
following assumptions:
\begin{enumerate}
\item$\hat{\mathbf{z}}=\arg\max_{z \in\mathbb{R}^{q}} [\log g(\mathbf
{y}|\mathbf{z}; \bolds\theta) + \log h(\mathbf{z})]$ exists for
all $l=1,\ldots, n$.
\item$\ell(\bolds\theta)$ is a well-defined function under these
regularity conditions:
\begin{enumerate}[$R_{4}$.]
\item[$R_{1}$.] $\ell(\bolds\theta)$ has a unique maximum at
$\bolds\theta_{0} \in\Theta$;
\item[$R_{2}$.] $\Theta$ is compact;
\item[$R_{3}$.] $\ell(\bolds\theta)$ is continuous;
\item[$R_{4}$.] the empirical approximated log-likelihood function
$\tilde{\ell}(\bolds\theta)$ converges uniformly in probability
to $\ell(\bolds\theta)$.
\end{enumerate}
\end{enumerate}
It has to be noticed that, under concavity of the objective
function $\tilde{\ell}(\bolds\theta)$, compactness ($R_{2}$) can
be replaced by the assumption that
\begin{enumerate}[$R_{2b}$.]
\item[$R_{2b}$.] the true parameter value $\bolds\theta_{0}$ is
an interior point of the parameter space, and the estimator $\hat
{\bolds\theta}$ is an interior point in a neighborhood containing
$\bolds\theta_{0}$ (see, e.g., Theorem 2.7 of Newey and
McFadden \cite{NewMcF94}).
\end{enumerate}
%
%
Let $\tilde{S}(\cdot)$ denote the approximated score vector
according to the approximations (\ref{sco2}); then we obtain
%
%
\begin{eqnarray}\label{appr}
&&\sum_{l=1}^{n} E_{\mathbf{z}|\mathbf
{y}}
\bigl[S_{l}(\hat{\bolds\theta}; \mathbf{z}_{l})\bigr]=
S(\hat{\bolds\theta})=\sum_{l=1}^{n}
\bigl\{S_{l}(\hat{\bolds\theta},\hat{\mathbf{z}}_{l})+
\cdots+ \RMO \bigl(p^{- [{k}/{3}+1 ]} \bigr) \bigr\}
\nonumber\\[-8pt]\\[-8pt]
&&\quad\Rightarrow \quad n^{-1}S(\hat{\bolds\theta})=n^{-1}
\tilde{S}(\hat{\bolds\theta}) + \RMO \bigl(p^{- [{k}/{3}+1 ]}
\bigr)\nonumber
\end{eqnarray}
since $\hat{\bolds\theta}$ is chosen such that $\sum_{l=1}^{n}
\tilde{E}_{\mathbf{z}|\mathbf{y}}[S_{l}(\hat{\bolds\theta};
\mathbf{z}_{l})]=\tilde{S}(\hat{\bolds\theta})=0$. Under the
regularity conditions in assumption 2 and provided that $(\hat
{\bolds\theta}-\bolds\theta_{0})=\RMo_{p}(1)$, we can apply a
Taylor series expansion in $S(\bolds\theta)$ around the true
parameter vector $\bolds\theta_{0}$:
%
%
\begin{equation}
\label{taysco}S(\hat{\bolds\theta})=S(\bolds\theta_{0})+ H
\bigl(\bolds\theta^{*}\bigr) (\hat{\bolds\theta}- \bolds
\theta_{0}),
\end{equation}
where $\bolds\theta^{*}$ lies on the segment joining $\bolds
\theta_{0}$ and $\hat{\bolds\theta}$, and
\[
H(\bolds\theta
^{*})=\frac{\partial S(\bolds\theta)}{\partial\theta
}\bigg|_{\bolds\theta=\bolds\theta^{*}}= \sum_{l=1}^{n}\frac
{\partial S_{l}(\bolds\theta, \hat{\mathbf
{z}}_{l})}{\partial\theta}\bigg|_{\bolds\theta=\bolds
\theta^{*}}=\sum_{l=1}^{n}H_{l}\bigl(\bolds\theta^{*}, \hat{\mathbf
{z}}_{l}\bigr).
\]

From equations (\ref{appr}) and (\ref{taysco}), we obtain
\begin{eqnarray*}
&&
(\hat{\bolds\theta} - \bolds\theta_{0})=-
\Biggl\{n^{-1}\sum_{l=1}^{n}H_{l}
\bigl(\bolds\theta^{*}, \hat{\mathbf{z}}_{l}\bigr)
\Biggr\}^{-1} \bigl\{n^{-1} \bigl[S(\bolds\theta
_{0})-S(\hat{\bolds\theta}) \bigr] \bigr\}
\\
&&\quad\Rightarrow\quad(\hat{\bolds\theta} - \bolds\theta_{0})=
-\Biggl\{n^{-1}\sum_{l=1}^{n}H_{l}
\bigl(\bolds\theta^{*}, \hat{\mathbf{z}}_{l}\bigr)
\Biggr\} ^{-1} \bigl[n^{-1}S(\bolds\theta_{0})+\RMO
\bigl(p^{- [{k}/{3}+1 ]} \bigr) \bigr].
\end{eqnarray*}
In addition, under assumption 2, we have that, as $n \rightarrow\infty
$, $n^{-1}H(\bolds\theta^{*})\rightarrow^{p}E_{\mathbf
{y}}[H(\bolds\theta_{0})]$, where the expectation is taken with
respect to $f(\mathbf{y}; \bolds\theta)$, and $H(\bolds
\theta^{*})=\sum_{l=1}^{n}H_{l}(\bolds\theta^{*}, \hat{\mathbf
{z}}_{l})$. By further assuming that $E_{\mathbf{y}} \{
H(\bolds\theta_{0}) \}$ is non-singular, we obtain
\[
\bigl\{n^{-1}H\bigl(\bolds\theta^{*}\bigr) \bigr
\}^{-1}\rightarrow^{p}E_{\mathbf{y}} \bigl\{H(\bolds
\theta_{0}) \bigr\}^{-1}.
\]
It follows that
\begin{eqnarray*}
(\hat{\bolds\theta}
- \bolds\theta_{0})&=&-E_{\mathbf
{y}} \bigl[H(\bolds
\theta_{0}) \bigr]^{-1} \bigl[n^{-1}S(\bolds
\theta_{0})+\RMO \bigl(p^{- [{k}/{3}+1
]} \bigr) \bigr]
\\
&=&\RMO_{p} \bigl[\max\bigl(n^{-1/2},p^{- [{k}/{3}+1 ]}
\bigr) \bigr],
\end{eqnarray*}
where in the last step we use the fact that, under the regularity
conditions 1, $n^{-1}S(\bolds\theta_{0})=\RMO_{p}(n^{-1/2})$, and
$E_{\mathbf{y}} \{H(\bolds\theta_{0}) \}=\RMO_{p}(1)$.

\section{Development of the adaptive ML estimators for binary manifest
variables}\label{app3}

Let $\mathbf{y}=(y_{1},\ldots, y_{p})^{T}$ be a vector of observed
binary variables, having a Bernoulli distribution with expectation $\pi
_{j}(\mathbf{z}), j=1,\ldots,p$. Using the canonical link function for
Bernoulli distribution, we have
\[
\pi_{j}(\mathbf{z})=\frac{\exp(\alpha_{0j}+ \bolds\alpha
_{j}^{T}\mathbf{z} )}{1+\exp(\alpha_{0j}+ \bolds\alpha
_{j}^{T}\mathbf{z} )}.
\]
The scale parameter $\phi_{j}=1$, such that the conditional
distribution of each observed binary item given the latent variables
$\mathbf{z}$ is
\[
g_{j}(y_{j}|\mathbf{z}; \bolds\theta)=\exp
\bigl[y_{j} \bigl(\alpha_{0j}+ \bolds
\alpha_{j}^{T}\mathbf{z} \bigr)-\log\bigl(1+\exp\bigl(
\alpha_{0j}+ \bolds\alpha_{j}^{T}\mathbf{z}
\bigr) \bigr) \bigr],\qquad j=1,\ldots, p.
\]
It follows that the approximated log-likelihood function (\ref
{approlike}) results
%
%
\begin{eqnarray}\label{empapplike}
\tilde{\ell}(\bolds\theta)&=& \sum_{l=1}^{n}
\log\Biggl[2^{q/2}|\mathbf{T}_{l}| \nonumber\\
&&\hspace*{32.7pt}{}\times\sum
_{t_{1},\ldots, t_{q}} \exp\Biggl(\sum_{i=1}^{p}y_{jl}
\bigl(\alpha_{0j}+\bolds\alpha_{j}^{T}\mathbf
{z}_{l, t_{1},\ldots,t_{q}}^{*} \bigr)\nonumber\\
&&\hspace*{89.4pt}{}-\sum_{i=1}^{p}
\log\bigl(1+\exp\bigl(\alpha_{0j}+ \bolds\alpha_{j}^{T}
\mathbf{z}_{l, t_{1},\ldots,t_{q}}^{*} \bigr) \bigr) \Biggr)
\nonumber\\
&&\hspace*{65.5pt}{}\times(2\pi)^{-q/2}\exp\biggl(-\frac{1}{2}
\mathbf{z}_{l, t_{1},\ldots,t_{q}}^{*T}\mathbf{z}_{l, t_{1},\ldots,t_{q}}^{*}
\biggr)w_{t_{1}}^{*}\cdots w_{t_{q}}^{*}
\Biggr]
\\
&=& \sum_{l=1}^{n} \Biggl\{-
\frac{q}{2} \log\pi\nonumber\\
&&\hspace*{19.5pt}{}+ \log|\mathbf{T}_{l}| +\log\Biggl[ \sum
_{t_{1},\ldots, t_{q}} \exp\Biggl(\sum_{j=1}^{p}y_{jl}
\bigl(\alpha_{0j}+\bolds\alpha_{j}^{T}\mathbf
{z}_{l, t_{1},\ldots,t_{q}}^{*} \bigr)
\nonumber\\
&&\hspace*{135.8pt}{}-\sum_{j=1}^{p}\log
\bigl(1+\exp\bigl(\alpha_{0j}+ \bolds\alpha_{j}^{T}
\mathbf{z}_{l, t_{1},\ldots,t_{q}}^{*} \bigr) \bigr)\nonumber\\
&&\hspace*{135.8pt}{}-\frac{1}{2}
\mathbf{z}_{l, t_{1},\ldots,t_{q}}^{*T}\mathbf{z}_{l, t_{1},\ldots,t_{q}}^{*}
\Biggr) w_{t_{1}}^{*}\cdots w_{t_{q}}^{*}
\Biggr]\Biggr\},\nonumber
\end{eqnarray}
where the AGH nodes and weights are derived by the classical
Gauss--Hermite nodes $z_{t_{k}}$ and weights $w_{t_{k}}, k=1,\ldots, q$,
as follows
\[
\mathbf{z}^{*}_{l, t_{1},\ldots, t_{q}}=\bigl(z_{l, t_{1}}^{*},\ldots,
z_{l, t_{q}}^{*}\bigr)^{T}=
\sqrt{2}\mathbf{T}_{l}(z_{t_{1}},\ldots, z_{t_{q}})^{T}+
\hat{\mathbf{z}}_{l}
\]
and
\[
w_{t_{k}}^{*}=w_{t_{k}}\exp\bigl[z_{t_{k}}^{2}
\bigr]
\]
with $\mathbf{T}_{l}$ derived by the Cholesky factorization of the
matrix $\bolds\Psi_{l}$, that is, $\bolds\Psi_{l}=\mathbf
{T}_{l}\mathbf{T}_{l}^{T}$. The modes $\hat{\mathbf{z}}_{l}$ are
obtained for each subject through the iterative scheme
\[
\hat{\mathbf{z}}_{l}^{\mathrm{it}+1}=\hat{\mathbf{z}}_{l}^{\mathrm{it}}+
\bolds\Psi_{l}^{\mathrm{it}}L\bigl(\hat{\mathbf{z}}_{l}^{\mathrm{it}}
\bigr),
\]
where ``it'' denotes the iteration counter,
\[
L\bigl(\hat{\mathbf{z}}_{l}^{\mathrm{it}}\bigr)=-
\frac{\partial[\log g(\mathbf
{y}_{l} \mid\mathbf{z}_{l}; \bolds\theta)+ \log h(\mathbf
{z}_{l}) ]}{\partial\mathbf{z}_{l}^{T}}\bigg|_{\mathbf
{z}_{l}=\hat{\mathbf{z}}_{l}^{\mathrm{it}}}=-\sum_{j=1}^{p}
\bolds\alpha_{j} \biggl[y_{jl}-\frac{\exp(\alpha_{0j}+
\bolds\alpha
_{j}^{T}\hat{\mathbf{z}}_{l}^{\mathrm{it}} )}{1+\exp(\alpha_{0j}+
\bolds\alpha_{j}^{T}\hat{\mathbf{z}}_{l}^{\mathrm{it}} )}
\biggr]+\hat{\mathbf{z}}_{l}^{\mathrm{it}}
\]
and
\[
\bolds\Psi_{l}^{-1}=-\frac{\partial^{2} [\log g(\mathbf
{y}_{l} \mid\mathbf{z}_{l}; \bolds\theta)+ \log h(\mathbf
{z}_{l}) ]}{\partial\mathbf{z}_{l}^{T}\partial\mathbf
{z}_{l}}\bigg|_{\mathbf{z}_{l}=\hat{\mathbf{z}}_{l}^{\mathrm{it}}}=
\sum_{j=1}^{p}\bolds
\alpha_{j}\bolds\alpha_{j}^{T}
\frac{\exp
(\alpha_{0j}+ \bolds\alpha_{j}^{T}\hat{\mathbf
{z}}_{l}^{\mathrm{it}} )}{1+\exp(\alpha_{0j}+ \bolds\alpha
_{j}^{T}\hat{\mathbf{z}}_{l}^{\mathrm{it}} )}+\mathbf{I}.
\]
%
\subsection{Regularity conditions for adaptive $M$-estimators in presence
of binary data}\label{secA.1}
Since the general theory of the $M$-estimators is here applied to a
particular family of GLLVM, the regularity conditions on the
log-likelihood function $\ell(\bolds\theta)$ given in Appendix \ref{app2}
should be checked for the particular distribution of each observed variable
(Huber, Ronchetti and Victoria-Feser \cite{HubRonFe04}). For classical
Laplace-based estimators, a formal proof of these conditions in
presence of ordinal manifest variables is given by Huber, Scaillet and
Victoria-Feser \cite{HubScaFe09}. Following the main lines of that paper,
we now prove how the empirical approximated log-likelihood (\ref
{empapplike}) satisfies the regularity conditions for consistency and
asymptotic normality of the corresponding $M$-estimators.

At this
regard, we make use of the Lemma 2.2 by Newey and McFadden \cite
{NewMcF94}, according to which $\ell(\bolds\theta)$ has a unique
maximum at $\bolds\theta_{0} \in\Theta$ (condition $R_{1}$) if:
\begin{enumerate}[$a_{2}$.]
\item[$a_{1}$.] $\bolds\theta_{0}$ is identified, that is, if
$\bolds\theta\neq\bolds\theta_{0}$, $\bolds\theta
\in\Theta$, then $\ell(\bolds\theta) \neq\ell(\bolds\theta
_{0})$, and

\item[$a_{2}$.] $E [|\tilde{\ell}(\bolds\theta)| ] <
\infty$.
\end{enumerate}
[$a_{1}$] Under our assumptions for the latent variables, $\bolds
\theta_{0}$ is identified.

\noindent[$a_{2}$] Let $\mathbf{z}^{*}$ be $\mathbf{z}_{l,
t_{1},\ldots,t_{q}}^{*}$ and let $K(\mathbf{z}^{*})$ denote $\log
(1+\exp(\alpha_{0j}+ \bolds\alpha_{j}^{T}\mathbf{z}^{*} ) )$. We
recall that $|\log(x)| \leq k (|x|+1 )$ for a constant $k \geq3$ and
for any $x >0$, and that $|\exp(x)|\leq\exp(|x|)$ for any $x
\in\mathbb{R}$. Hence, based on (\ref{empapplike}),
%
%
\begin{eqnarray}
\label{expon}
&&\sum_{l=1}^{n}\Biggl\llvert\log
\Biggl[ \sum_{t_{1},\ldots, t_{q}} \exp\Biggl(\sum
_{j=1}^{p}y_{jl} \bigl(\alpha_{0j}+
\bolds\alpha_{j}^{T}\mathbf{z}^{*} \bigr)-
\sum_{j=1}^{p}K \bigl(\mathbf{z}^{*}
\bigr)-\frac
{1}{2}\mathbf{z}^{*T}\mathbf{z}^{*}
\Biggr)w_{t_{1}}^{*}\cdots w_{t_{q}}^{*}
\Biggr]\Biggr\rrvert
\nonumber\\[-8pt]\\[-8pt]
&&\quad\leq \sum_{l=1}^{n} k \Biggl(
\sum_{t_{1},\ldots, t_{q}} \exp\Biggl(\sum
_{j=1}^{p} \bigl\llvert y_{jl} \bigl(
\alpha_{0j}+\bolds\alpha_{j}^{T}
\mathbf{z}^{*} \bigr) \bigr\rrvert+\sum_{j=1}^{p}
\bigl\llvert K \bigl(\mathbf{z}^{*} \bigr)\bigr\rrvert+
\frac{1}{2}\bigl\llVert\mathbf{z}^{*}\bigr\rrVert
\Biggr)w_{t_{1}}^{*}\cdots w_{t_{q}}^{*} +1
\Biggr).\qquad\nonumber
\end{eqnarray}
It can be noticed that $|K(\mathbf{z}^{*})| \leq\log2$ if $
(\alpha_{0j}+\bolds\alpha_{j}^{T}\mathbf{z}^{*} ) < 0$, and
$|K(\mathbf{z}^{*})| \leq|k_{1}||\bolds\alpha_{j}|\|\mathbf
{z}^{*}\|+ |k_{1}||\alpha_{0j}|+|k_{2}|= |k_{1}||\bolds\alpha
_{j}|\|\mathbf{z}^{*}\|+ \mbox{Const}$ if $ (\alpha_{0j}+\bolds\alpha
_{j}^{T}\mathbf{z}^{*} ) > 0$ (using $\log(1+x) \leq k_{1}\log
(x)+k_{2}$ for any constant $k_{1} \geq\frac{1}{2}$ and $k_{2} > 1$).
Furthermore, $|y_{jl} (\alpha_{0j}+\bolds\alpha_{j}^{T}\mathbf
{z}^{*} )| \leq|y_{jl}||\bolds\alpha_{j}| \|\mathbf
{z}^{*}\| + |y_{jl}||\alpha_{0j}|= |y_{jl}||\bolds\alpha_{j}|
\|\mathbf{z}^{*}\| + \mbox{Const}$. Using the definition of $\mathbf{z}^{*}$,
we deduce that $E \|\mathbf{z}^{*}\|= \sum_{t_{1},\ldots, t_{q}}\mathbf
{z}_{l, t_{1},\ldots,t_{q}}^{*T}\mathbf{z}_{l, t_{1},\ldots,t_{q}}^{*}
w_{t_{1}}^{*}\cdots w_{t_{q}}^{*} < \infty$. Hence, based on log-normal
moments, (\ref{expon}) is finite. Besides 
$|{\log\det(\mathbf{T}_{l}) }| < \mbox{Const}$, such that $E [|\tilde{\ell
}(\bolds\theta)| ]<\infty$.

Since the data are i.i.d. and
$\Theta$ is compact (condition $R_{2}$), $\tilde{\ell}(\bolds
\theta)$ is continuous at each $\bolds\theta$ with probability
one, and there is a function of the latent variables $d(\mathbf
{z}^{*})$ with $|\tilde{\ell}(\bolds\theta)| \leq d(\mathbf
{z}^{*})$ such that $E [d(\mathbf{z}^{*}) ]< \infty$ (cf.
proof of condition $a_{2}$). So, we deduce that $E [\tilde{\ell
}(\bolds\theta) ]=\ell(\bolds\theta)$ is continuous
(condition $R_{3}$) and that $\tilde{\ell}(\bolds\theta)$
converges uniformly in probability to that quantity (condition $R_{4}$)
(see Lemma 2.4 by Newey and McFadden \cite{NewMcF94}).

Asymptotic normality of the estimators imposes conditions on the
Hessian of the empirical approximated log-likelihood function, that
should be verified. Based on (\ref{empapplike}), by the computing
the explicit expression of $\frac{\partial^{2} \tilde{\ell}(\bolds
\theta)}{\partial\bolds\theta\,\partial\bolds\theta^{T}}$,
it can be easily shown that there is a function $d(\mathbf{z}^{*})$
with $\llvert\frac{\partial^{2} \tilde{\ell}(\bolds\theta
)}{\partial\bolds\theta\,\partial\bolds\theta^{T}}\rrvert<
d(\mathbf{z}^{*})$ such that $E [d(\mathbf{z}^{*}) ]< \infty$
(as in the proof of condition $a_{2}$). As before, making use of Lemma
2.4 by Newey and McFadden \cite{NewMcF94}, since the data are i.i.d.
and $\Theta$ is compact, $\frac{\partial^{2} \tilde{\ell}(\bolds
\theta)}{\partial\bolds\theta\,\partial\bolds\theta^{T}}$
is continuous\vspace*{1pt} at each $\bolds\theta$ with probability one, and we
can deduce that $E [\frac{\partial^{2} \tilde{\ell}(\bolds
\theta)}{\partial\bolds\theta\,\partial\bolds\theta
^{T}} ]$ is continuous and the Hessian of the empirical
approximated log-likelihood converges uniformly in probability to that quantity.


\subsection{Score functions and second order Laplace estimators}\label{secA.2}
In this specific case, the complete data score functions (\ref{sco1})
are given by
\begin{eqnarray*}
S_{l}\bigl(\alpha_{0j}; \mathbf{z}^{*}_{l, t_{1},\ldots, t_{q}}
\bigr)&=& \biggl[y_{jl}-\frac{\exp(\alpha_{0j}+ \bolds\alpha
_{j}^{T}\mathbf
{z}^{*}_{l, t_{1},\ldots, t_{q}} )}{1+\exp(\alpha_{0j}+
\bolds\alpha_{j}^{T}\mathbf{z}^{*}_{l, t_{1},\ldots, t_{q}}
)} \biggr],
\\
S_{l}\bigl(\bolds\alpha_{j}; \mathbf{z}^{*}_{l, t_{1},\ldots,
t_{q}}
\bigr)&=& \mathbf{z}^{*}_{l, t_{1},\ldots, t_{q}} \biggl[y_{jl}-
\frac{\exp
(\alpha_{0j}+ \bolds\alpha_{j}^{T}\mathbf{z}^{*}_{l, t_{1},\ldots,
t_{q}} )}{1+\exp(\alpha_{0j}+ \bolds\alpha
_{j}^{T}\mathbf{z}^{*}_{l, t_{1},\ldots, t_{q}} )} \biggr].
\end{eqnarray*}
The corresponding score equations have not closed form solutions, and a
quasi-Newton procedure is used to solve implicit equations.

In the simulation study, the performance of the adaptive-based
estimators has been compared with second order Laplace estimators.
According to (\ref{Ei2}), the latter have been derived by
maximizing the following approximated log-likelihood function
\[
\tilde{\ell}(\bolds\theta)=\sum_{l=1}^{n}
\log\bigl\{(2\pi)^{q/2}|\bolds\Psi_{l}^{1/2}
\exp\bigl[-L(\hat{\mathbf{z}}_{l})\bigr] \bigl[1+c_{1}p^{-1}+\RMO
\bigl(p^{-2}\bigr) \bigr] \bigr\},
\]
where the individual modes $\hat{\mathbf{z}}_{l}$ are obtained through
the iterative scheme defined above, and
\[
c_{1}=\sum_{m=2}^{3}
\frac{(-1)^{m-1}}{(2m)!}L_{p_{1}}(\hat{\mathbf{z}}_{l}) \cdots
L_{p_{m-1}}(\hat{\mathbf{z}}_{l})L^{q_{1}}(\hat{\mathbf
{z}}_{l}) \cdots L^{q_{m}}(\hat{\mathbf{z}}_{l})
\]
with $p_{1}|\cdots|p_{m-1}$ be a partition of $2m$ indices into $m-1$
blocks, each of size 3 or more, and $q_{1}|\cdots|q_{m}$ is a
partition of $2m$ indices into $m$ blocks, each of size 2. In
particular, following the notation by Raudenbush, Yang and Yosef \cite{Rade00},
\[
c_{1}=-\tfrac{1}{8}\operatorname{vec}^{T}[\bolds
\Psi_{l} \otimes\bolds\Psi_{l}] \operatorname{vec}
\bigl[L^{(4)}(\hat{\mathbf{z}}_{l})\bigr] +
\tfrac
{5}{24}\operatorname{vec}^{T}[\bolds\Psi_{l} \otimes
\bolds\Psi_{l} \otimes\bolds\Psi_{l}] \operatorname{vec}
\bigl[L^{(3)}(\hat{\mathbf{z}}_{l}) \otimes L^{(3)}(
\hat{\mathbf{z}}_{l})\bigr],
\]
where
\[
L^{(3)}(\hat{\mathbf{z}}_{l})=-\sum
_{j=1}^{p}\operatorname{vec}\bigl(\bolds\alpha_{j}
\bolds\alpha_{j}^{T}\bigr)\bolds
\alpha_{j}^{T}\frac{\exp
(\alpha_{0j}+ \bolds\alpha_{j}^{T}\hat{\mathbf{z}}_{l} )
[1-\exp(\alpha_{0j}+ \bolds\alpha_{j}^{T}\hat{\mathbf
{z}}_{l} ) ]}{ [1+\exp(\alpha_{0j}+ \bolds
\alpha_{j}^{T}\hat{\mathbf{z}}_{l} ) ]^{3}}
\]
and
\begin{eqnarray*}
L^{(4)}(\hat{\mathbf{z}}_{l})&=&-\sum
_{j=1}^{p}\operatorname{vec}\bigl[\operatorname{vec}\bigl(\bolds\alpha
_{j}\bolds\alpha_{j}^{T}\bigr)\bolds
\alpha_{j}^{T}\bigr]\\
&&\hspace*{22.4pt}{}\times\bolds\alpha_{j}^{T}
\frac{\exp(\alpha_{0j}+ \bolds\alpha
_{j}^{T}\hat{\mathbf{z}}_{l} ) [1-4\exp(\alpha_{0j}+
\bolds\alpha_{j}^{T}\hat{\mathbf{z}}_{l} )+\exp(\alpha
_{0j}+ \bolds\alpha_{j}^{T}\hat{\mathbf{z}}_{l} )^{2}
]}{ [1+\exp(\alpha_{0j}+ \bolds\alpha_{j}^{T}\hat{\mathbf
{z}}_{l} ) ]^{4}}.
\end{eqnarray*}
As for the adaptive-based estimators, the score equations of both the
intercepts and factor loadings have not closed form solutions, and a
quasi-Newton procedure has been used to solve implicit equations.
\end{appendix}



%

\printhistory


\begin{thebibliography}{35}

\bibitem{BarCox89}
%
\begin{bbook}[mr]
\bauthor{\bsnm{Barndorff-Nielsen},~\bfnm{O.~E.}\binits{O.E.}} \AND
\bauthor{\bsnm{Cox},~\bfnm{D.~R.}\binits{D.R.}}
(\byear{1989}).
\btitle{Asymptotic Techniques for Use in Statistics}.
\bseries{Monographs on Statistics and Applied Probability}.
\blocation{London}: \bpublisher{Chapman \& Hall}.
\bid{mr={1010226}}
\bptok{imsref}%
\end{bbook}
%
\endbibitem


\bibitem{BiaCag12}
%
\begin{barticle}[mr]
\bauthor{\bsnm{Bianconcini},~\bfnm{Silvia}\binits{S.}} \AND
\bauthor{\bsnm{Cagnone},~\bfnm{Silvia}\binits{S.}}
(\byear{2012}).
\btitle{Estimation of generalized linear latent variable models via fully
exponential {L}aplace approximation}.
\bjournal{J. Multivariate Anal.}
\bvolume{112}
\bpages{183--193}.
\bid{doi={10.1016/j.jmva.2012.06.005}, issn={0047-259X}, mr={2957295}}
\bptok{imsref}%
\end{barticle}
%
\endbibitem


\bibitem{DeBr81}
%
\begin{bbook}[mr]
\bauthor{\bparticle{de} \bsnm{Bruijn},~\bfnm{N.~G.}\binits{N.G.}}
(\byear{1981}).
\btitle{Asymptotic Methods in Analysis},
\bedition{3rd} ed.
\blocation{New York}: \bpublisher{Dover}.
\bid{mr={0671583}}
\bptok{imsref}%
\end{bbook}
%
\endbibitem

\bibitem{Eva11}
%
\begin{barticle}[mr]
\bauthor{\bsnm{Evangelou},~\bfnm{Evangelos}\binits{E.}},
\bauthor{\bsnm{Zhu},~\bfnm{Zhengyuan}\binits{Z.}} \AND
\bauthor{\bsnm{Smith},~\bfnm{Richard~L.}\binits{R.L.}}
(\byear{2011}).
\btitle{Estimation and prediction for spatial generalized linear mixed models
using high order {L}aplace approximation}.
\bjournal{J. Statist. Plann. Inference}
\bvolume{141}
\bpages{3564--3577}.
\bid{doi={10.1016/j.jspi.2011.05.008}, issn={0378-3758}, mr={2817363}}
\bptok{imsref}%
\end{barticle}
%
\endbibitem

\bibitem{HubRonFe04}
%
\begin{barticle}[mr]
\bauthor{\bsnm{Huber},~\bfnm{Philippe}\binits{P.}},
\bauthor{\bsnm{Ronchetti},~\bfnm{Elvezio}\binits{E.}} \AND
\bauthor{\bsnm{Victoria-Feser},~\bfnm{Maria-Pia}\binits{M.P.}}
(\byear{2004}).
\btitle{Estimation of generalized linear latent variable models}.
\bjournal{J. R. Stat. Soc. Ser. B Stat. Methodol.}
\bvolume{66}
\bpages{893--908}.
\bid{doi={10.1111/j.1467-9868.2004.05627.x}, issn={1369-7412}, mr={2102471}}
\bptok{imsref}%
\end{barticle}
%
\endbibitem

\bibitem{HubScaFe09}
%
\begin{barticle}[mr]
\bauthor{\bsnm{Huber},~\bfnm{Philippe}\binits{P.}},
\bauthor{\bsnm{Scaillet},~\bfnm{Olivier}\binits{O.}} \AND
\bauthor{\bsnm{Victoria-Feser},~\bfnm{Maria-Pia}\binits{M.P.}}
(\byear{2009}).
\btitle{Assessing multivariate predictors of financial market
movements: A
latent factor framework for ordinal data}.
\bjournal{Ann. Appl. Stat.}
\bvolume{3}
\bpages{249--271}.
\bid{doi={10.1214/08-AOAS213}, issn={1932-6157}, mr={2668707}}
\bptok{imsref}%
\end{barticle}
%
\endbibitem


\bibitem{Joe08}
%
\begin{barticle}[mr]
\bauthor{\bsnm{Joe},~\bfnm{Harry}\binits{H.}}
(\byear{2008}).
\btitle{Accuracy of {L}aplace approximation for discrete response mixed
models}.
\bjournal{Comput. Statist. Data Anal.}
\bvolume{52}
\bpages{5066--5074}.
\bid{doi={10.1016/j.csda.2008.05.002}, issn={0167-9473}, mr={2526575}}
\bptok{imsref}%
\end{barticle}
%
\endbibitem

\bibitem{LeSp01}
%
\begin{barticle}[mr]
\bauthor{\bsnm{Lesaffre},~\bfnm{Emmanuel}\binits{E.}} \AND
\bauthor{\bsnm{Spiessens},~\bfnm{Bart}\binits{B.}}
(\byear{2001}).
\btitle{On the effect of the number of quadrature points in a logistic
random-effects model: An example}.
\bjournal{J. R. Stat. Soc. Ser. C. Appl. Stat.}
\bvolume{50}
\bpages{325--335}.
\bid{doi={10.1111/1467-9876.00237}, issn={0035-9254}, mr={1856328}}
\bptok{imsref}%
\end{barticle}
%
\endbibitem

\bibitem{LiuPie94}
%
\begin{barticle}[mr]
\bauthor{\bsnm{Liu},~\bfnm{Qing}\binits{Q.}} \AND
\bauthor{\bsnm{Pierce},~\bfnm{Donald~A.}\binits{D.A.}}
(\byear{1994}).
\btitle{A note on {G}auss--{H}ermite quadrature}.
\bjournal{Biometrika}
\bvolume{81}
\bpages{624--629}.
\bid{issn={0006-3444}, mr={1311107}}
\bptok{imsref}%
\end{barticle}
%
\endbibitem


\bibitem{Rue13}
%
\begin{barticle}[auto:STB|2013/06/05|13:45:01]
\bauthor{\bsnm{Martins},~\bfnm{T.~G.}\binits{T.G.}},
\bauthor{\bsnm{Simpson},~\bfnm{D.}\binits{D.}},
\bauthor{\bsnm{Lindgren},~\bfnm{F.}\binits{F.}} \AND
\bauthor{\bsnm{Rue},~\bfnm{H.}\binits{H.}}
(\byear{2013}).
\btitle{Bayesian computing with INLA: New features}.
\bjournal{Comput.
Statist. Data Anal.}
\bvolume{67}
\bpages{68--83}.
\bptok{imsref}%
\end{barticle}
%
\endbibitem


\bibitem{MouKno00}
%
\begin{barticle}[mr]
\bauthor{\bsnm{Moustaki},~\bfnm{Irini}\binits{I.}} \AND
\bauthor{\bsnm{Knott},~\bfnm{Martin}\binits{M.}}
(\byear{2000}).
\btitle{Generalized latent trait models}.
\bjournal{Psychometrika}
\bvolume{65}
\bpages{391--411}.
\bid{doi={10.1007/BF02296153}, issn={0033-3123}, mr={1792703}}
\bptok{imsref}%
\end{barticle}
%
\endbibitem

\bibitem{MuMu10}
%
\begin{bbook}[auto:STB|2013/06/05|13:45:01]
\bauthor{\bsnm{Muth\'{e}n},~\bfnm{L.~K.}\binits{L.K.}} \AND
\bauthor{\bsnm{Muth\'{e}n},~\bfnm{B.~O.}\binits{B.O.}}
(\byear{2010}).
\btitle{Mplus User's Guide},
\bedition{6th} ed.
\blocation{Los Angeles, CA}: \bpublisher{Muth\'{e}n \& Muth\'{e}n}.
\bptok{imsref}%
\end{bbook}
%
\endbibitem

\bibitem{NaySmi82}
%
\begin{barticle}[mr]
\bauthor{\bsnm{Naylor},~\bfnm{J.~C.}\binits{J.C.}} \AND
\bauthor{\bsnm{Smith},~\bfnm{A.~F.~M.}\binits{A.F.M.}}
(\byear{1982}).
\btitle{Applications of a method for the efficient computation of posterior
distributions}.
\bjournal{J. R. Stat. Soc. Ser. C. Appl. Stat.}
\bvolume{31}
\bpages{214--225}.
\bid{doi={10.2307/2347995}, issn={0035-9254}, mr={0694917}}
\bptok{imsref}%
\end{barticle}
%
\endbibitem

\bibitem{NewMcF94}
%
\begin{bincollection}[mr]
\bauthor{\bsnm{Newey},~\bfnm{Whitney~K.}\binits{W.K.}} \AND
\bauthor{\bsnm{McFadden},~\bfnm{Daniel}\binits{D.}}
(\byear{1994}).
\btitle{Large sample estimation and hypothesis testing}.
In \bbooktitle{Handbook of Econometrics, {V}ol. {IV}}
(\beditor{R.F. Engle} and \beditor{D. McFadden}, eds.).
\bseries{Handbooks in Econom.}
\bvolume{2}
\bpages{2111--2245}.
\blocation{Amsterdam}: \bpublisher{North-Holland}.
\bid{mr={1315971}}
\bptok{imsref}%
\end{bincollection}
%
\endbibitem

\bibitem{PiBa95}
%
\begin{barticle}[auto:STB|2013/06/05|13:45:01]
\bauthor{\bsnm{Pinhero},~\bfnm{J.~C.}\binits{J.C.}} \AND
\bauthor{\bsnm{Bates},~\bfnm{D.~M.}\binits{D.M.}}
(\byear{1995}).
\btitle{Approximation to the log-likelihood function in the nonlinear mixed
effects model}.
\bjournal{J. Comput. Graph. Statist.}
\bvolume{4}
\bpages{12--35}.
\bptok{imsref}%
\end{barticle}
%
\endbibitem

\bibitem{R10}
%
\begin{bmisc}[auto:STB|2013/06/05|13:45:01]
\borganization{R Development Core Team}
(\byear{2010}).
\bhowpublished{R: A language and environment for statistical computing. R
Foundation for Statistical Computing, Vienna}.
\bptok{imsref}%
\end{bmisc}
%
\endbibitem

\bibitem{RaSk12}
%
\begin{bbook}[auto:STB|2013/06/05|13:45:01]
\bauthor{\bsnm{Rabe-Hesketh},~\bfnm{S.}\binits{S.}} \AND
\bauthor{\bsnm{Skrondal},~\bfnm{A.}\binits{A.}}
(\byear{2012}).
\btitle{Multilevel and Longitudinal Modeling Using Stata},
\bedition{3rd} ed.
\blocation{College Station, TX}: \bpublisher{Stata Press}.
\bptok{imsref}%
\end{bbook}
%
\endbibitem

\bibitem{HesSkPi05}
%
\begin{barticle}[mr]
\bauthor{\bsnm{Rabe-Hesketh},~\bfnm{Sophia}\binits{S.}},
\bauthor{\bsnm{Skrondal},~\bfnm{Anders}\binits{A.}} \AND
\bauthor{\bsnm{Pickles},~\bfnm{Andrew}\binits{A.}}
(\byear{2005}).
\btitle{Maximum likelihood estimation of limited and discrete dependent
variable models with nested random effects}.
\bjournal{J.~Econometrics}
\bvolume{128}
\bpages{301--323}.
\bid{doi={10.1016/j.jeconom.2004.08.017}, issn={0304-4076}, mr={2189555}}
\bptok{imsref}%
\end{barticle}
%
\endbibitem

\bibitem{Rade00}
%
\begin{barticle}[mr]
\bauthor{\bsnm{Raudenbush},~\bfnm{Stephen~W.}\binits{S.W.}},
\bauthor{\bsnm{Yang},~\bfnm{Meng-Li}\binits{M.L.}} \AND
\bauthor{\bsnm{Yosef},~\bfnm{Matheos}\binits{M.}}
(\byear{2000}).
\btitle{Maximum likelihood for generalized linear models with nested random
effects via high-order, multivariate {L}aplace approximation}.
\bjournal{J. Comput. Graph. Statist.}
\bvolume{9}
\bpages{141--157}.
\bid{doi={10.2307/1390617}, issn={1061-8600}, mr={1826278}}
\bptok{imsref}%
\end{barticle}
%
\endbibitem

\bibitem{RizVer09}
%
\begin{barticle}[mr]
\bauthor{\bsnm{Rizopoulos},~\bfnm{Dimitris}\binits{D.}},
\bauthor{\bsnm{Verbeke},~\bfnm{Geert}\binits{G.}} \AND
\bauthor{\bsnm{Lesaffre},~\bfnm{Emmanuel}\binits{E.}}
(\byear{2009}).
\btitle{Fully exponential {L}aplace approximations for the joint
modelling of
survival and longitudinal data}.
\bjournal{J. R. Stat. Soc. Ser. B Stat. Methodol.}
\bvolume{71}
\bpages{637--654}.
\bid{doi={10.1111/j.1467-9868.2008.00704.x}, issn={1369-7412}, mr={2749911}}
\bptok{imsref}%
\end{barticle}
%
\endbibitem

\bibitem{Rue05}
%
\begin{bbook}[mr]
\bauthor{\bsnm{Rue},~\bfnm{H{\aa}vard}\binits{H.}} \AND
\bauthor{\bsnm{Held},~\bfnm{Leonhard}\binits{L.}}
(\byear{2005}).
\btitle{Gaussian {M}arkov Random Fields:
Theory and Applications}.
\bseries{Monographs on Statistics and Applied Probability}
\bvolume{104}.
\blocation{Boca Raton, FL}: \bpublisher{Chapman \& Hall/CRC}.
\bid{doi={10.1201/9780203492024}, mr={2130347}}
\bptok{imsref}%
\end{bbook}
%
\endbibitem

\bibitem{Rue09}
%
\begin{barticle}[mr]
\bauthor{\bsnm{Rue},~\bfnm{H{\aa}vard}\binits{H.}},
\bauthor{\bsnm{Martino},~\bfnm{Sara}\binits{S.}} \AND
\bauthor{\bsnm{Chopin},~\bfnm{Nicolas}\binits{N.}}
(\byear{2009}).
\btitle{Approximate {B}ayesian inference for latent {G}aussian models
by using
integrated nested {L}aplace approximations}.
\bjournal{J. R. Stat. Soc. Ser. B Stat. Methodol.}
\bvolume{71}
\bpages{319--392}.
\bid{doi={10.1111/j.1467-9868.2008.00700.x}, issn={1369-7412}, mr={2649602}}
\bptok{imsref}%
\end{barticle}
%
\endbibitem

\bibitem{SchBo05}
%
\begin{barticle}[mr]
\bauthor{\bsnm{Schilling},~\bfnm{Stephen}\binits{S.}} \AND
\bauthor{\bsnm{Bock},~\bfnm{R.~Darrell}\binits{R.D.}}
(\byear{2005}).
\btitle{High-dimensional maximum marginal likelihood item factor
analysis by
adaptive quadrature}.
\bjournal{Psychometrika}
\bvolume{70}
\bpages{533--555}.
\bid{doi={10.1007/s11336-003-1141-x}, issn={0033-3123}, mr={2272503}}
\bptok{imsref}%
\end{barticle}
%
\endbibitem

\bibitem{Shu97}
%
\begin{barticle}[auto:STB|2013/06/05|13:45:01]
\bauthor{\bsnm{Shun},~\bfnm{Z.}\binits{Z.}}
(\byear{1997}).
\btitle{Another look at the salamander mating data: A modified Laplace
approximation approach}.
\bjournal{J. Amer. Statist. Assoc.}
\bvolume{92}
\bpages{341--349}.
\bptok{imsref}%
\end{barticle}
%
\endbibitem

\bibitem{ShuMc95}
%
\begin{barticle}[mr]
\bauthor{\bsnm{Shun},~\bfnm{Zhenming}\binits{Z.}} \AND
\bauthor{\bsnm{McCullagh},~\bfnm{Peter}\binits{P.}}
(\byear{1995}).
\btitle{Laplace approximation of high-dimensional integrals}.
\bjournal{J. R. Stat. Soc. Ser. B Stat. Methodol.}
\bvolume{57}
\bpages{749--760}.
\bid{issn={0035-9246}, mr={1354079}}
\bptok{imsref}%
\end{barticle}
%
\endbibitem

\bibitem{SkRa04}
%
\begin{bbook}[mr]
\bauthor{\bsnm{Skrondal},~\bfnm{Anders}\binits{A.}} \AND
\bauthor{\bsnm{Rabe-Hesketh},~\bfnm{Sophia}\binits{S.}}
(\byear{2004}).
\btitle{Generalized Latent Variable Modeling:
Multilevel, Longitudinal, and Structural Equation Models}.
\bseries{Interdisciplinary Statistics}.
\blocation{Boca Raton, FL}: \bpublisher{Chapman \& Hall/CRC}.
\bid{doi={10.1201/9780203489437}, mr={2059021}}
\bptok{imsref}%
\end{bbook}
%
\endbibitem

\bibitem{TauHus91}
%
\begin{barticle}[mr]
\bauthor{\bsnm{Tauchen},~\bfnm{George}\binits{G.}} \AND
\bauthor{\bsnm{Hussey},~\bfnm{Robert}\binits{R.}}
(\byear{1991}).
\btitle{Quadrature-based methods for obtaining approximate solutions to
nonlinear asset pricing models}.
\bjournal{Econometrica}
\bvolume{59}
\bpages{371--396}.
\bid{doi={10.2307/2938261}, issn={0012-9682}, mr={1097533}}
\bptok{imsref}%
\end{barticle}
%
\endbibitem

\bibitem{TieKa89}
%
\begin{barticle}[mr]
\bauthor{\bsnm{Tierney},~\bfnm{Luke}\binits{L.}},
\bauthor{\bsnm{Kass},~\bfnm{Robert~E.}\binits{R.E.}} \AND
\bauthor{\bsnm{Kadane},~\bfnm{Joseph~B.}\binits{J.B.}}
(\byear{1989}).
\btitle{Fully exponential {L}aplace approximations to expectations and
variances of nonpositive functions}.
\bjournal{J. Amer. Statist. Assoc.}
\bvolume{84}
\bpages{710--716}.
\bid{issn={0162-1459}, mr={1132586}}
\bptok{imsref}%
\end{barticle}
%
\endbibitem

\bibitem{Von96}
%
\begin{barticle}[mr]
\bauthor{\bsnm{Vonesh},~\bfnm{Edward~F.}\binits{E.F.}}
(\byear{1996}).
\btitle{A note on the use of {L}aplace's approximation for nonlinear
mixed-effects models}.
\bjournal{Biometrika}
\bvolume{83}
\bpages{447--452}.
\bid{doi={10.1093/biomet/83.2.447}, issn={0006-3444}, mr={1439795}}
\bptok{imsref}%
\end{barticle}
%
\endbibitem

\end{thebibliography}
\end{document}